\documentclass[11pt,a4paper]{article}
\usepackage[utf8]{inputenc}
\usepackage{graphicx,verbatim,array,multicol,courier}
\usepackage{amsmath,amssymb,amsthm,mathrsfs, mathtools}
\usepackage{color, graphics, graphicx, siunitx}
\usepackage{amsfonts, dsfont, bm}
\usepackage{natbib} 
\usepackage{geometry}
\usepackage{enumitem,float}
\usepackage{lscape}
\usepackage[pdftex,bookmarks,colorlinks]{hyperref}
\usepackage[dvipsnames]{xcolor}
\usepackage{tabularx}
\usepackage{siunitx}
\usepackage[labelformat=simple]{subfig}
\usepackage[linesnumbered,ruled,vlined]{algorithm2e}
\usepackage{booktabs}

\newlist{inparaenum}{enumerate}{2}
\setlist[inparaenum,1]{label=(\alph*)}
\setlist[inparaenum,2]{label=(\roman{inparaenumi}\emph{\alph*})}

\makeatletter
\def\adl@drawiv#1#2#3{%
        \hskip.5\tabcolsep
        \xleaders#3{#2.5\@tempdimb #1{1}#2.5\@tempdimb}%
                #2\z@ plus1fil minus1fil\relax
        \hskip.5\tabcolsep}
\newcommand{\cdashlinelr}[1]{%
  \noalign{\vskip\aboverulesep
           \global\let\@dashdrawstore\adl@draw
           \global\let\adl@draw\adl@drawiv}
  \cdashline{#1}
  \noalign{\global\let\adl@draw\@dashdrawstore
           \vskip\belowrulesep}}
\makeatother

\theoremstyle{plain}
\newtheorem{axiom}{Axiom}

\newtheorem{cond}[axiom]{Condition}
\newtheorem{theorem}{Theorem}[section]

\newtheorem{proposition}[theorem]{Proposition}
\newtheorem{corollary}[theorem]{Corollary}
\theoremstyle{definition}

\newtheorem*{example}{Example}

\newtheorem{rem}[theorem]{Remark}

\newcommand{\diff}{\mathrm{d}}

\definecolor{navy}{rgb}{0,0,0.502}
\definecolor{brown}{rgb}{0.59, 0.29, 0.0}
\def\indic{\mathds{1}}

\hypersetup{colorlinks,%
	citecolor=blue,%
	filecolor=green,%
	linkcolor=red,%
	urlcolor=violet,%
}

\def\indic{\mathds{1}}

\newcommand{\lik}{{\mathscr{L}}}
\newcommand{\Pmu}{{\mathcal{P}}}
\newcommand{\Psp}{{\mathscr{P}}}

\newcommand{\Supp}{{\mathcal{S}}}
\newcommand{\bbP}{\mathbb{P}}

\newcommand{\Real}{\mathbb{R}}

\newcommand{\DoA}{\mathcal{D}}
\newcommand{\Prob}{\mathbb{P}}
\newcommand{\Expect}{\mathbb{E}}

\newcommand{\bftheta}{{\boldsymbol{\theta}}}
\newcommand{\bfvartheta}{{\boldsymbol{\vartheta}}}

\newcommand{\bfmu}{{\boldsymbol{\mu}}}

\newcommand{\bfSigma}{{\boldsymbol{\Sigma}}}
\newcommand{\bfDelta}{{\boldsymbol{\Delta}}}

\newcommand{\bfJ}{{\boldsymbol{J}}}

\newcommand{\bfX}{{\boldsymbol{X}}}
\newcommand{\bfY}{{\boldsymbol{Y}}}

\newcommand{\bfb}{{\boldsymbol{b}}}
\newcommand{\bfI}{{\boldsymbol{I}}}
\newcommand{\bfy}{{\boldsymbol{y}}}
\newcommand{\bfx}{{\boldsymbol{x}}}

\newcommand{\bfv}{{\boldsymbol{v}}}

\newcommand{\bfw}{{\boldsymbol{w}}}

\title{Asymptotic theory for Bayesian inference and prediction: from the ordinary to a conditional Peaks-Over-Threshold method}
\author{Cl\'ement Dombry, Simone A. Padoan and Stefano Rizzelli}

\begin{document}

\maketitle
\begin{abstract}
The Peaks Over Threshold (POT) method is the most popular statistical method for the analysis of univariate extremes. Even though there is a rich applied literature on Bayesian inference for the POT, the asymptotic theory for such proposals is missing. Even more importantly, the ambitious and challenging problem of predicting future extreme events according to a proper predictive statistical approach has received no attention to date. In this paper we fill this gap by developing the asymptotic theory of posterior distributions (consistency, contraction rates, asymptotic normality and asymptotic coverage of credible intervals) and prediction within the Bayesian framework in the POT context. We extend this asymptotic theory to account for cases where the focus is on the tail properties of the conditional distribution of a response variable given a vector of random covariates. To enable accurate predictions of extreme events more severe than those previously observed, we derive the posterior predictive distribution as an estimator of the conditional distribution of an out-of-sample random variable, given that it exceeds a sufficiently high threshold. We establish Wasserstein consistency of the posterior predictive distribution under both the unconditional and covariate-conditional approaches and derive its contraction rates. Simulations show the good performances of the proposed Bayesian inferential methods. The analysis of the change in the frequency of financial crises over time shows the utility of our methodology.
\end{abstract}

%
%
%
\section{Introduction}\label{sec:introduction}
%
%

%
%
\subsection{Statistical model and its challenges}\label{sec:stat_models}
%
%
%
The mission of Extreme Values Theory (EVT) is modelling and predicting  future events that are much more exceptional than those experienced in the past. Accomplishing this goal is undoubtedly important in many applied fields and for this purpose EVT develops tools for supporting risk assessment.
Here we focus on the Peaks Over Threshold (POT) method which is arguably the most popular approach in the univariate case. 
In the first part, we work with a random variable $Y$ with a generic distribution $F$. 
Under weak conditions the distribution of $Y-t\mid Y>t$, for a high threshold $t$ that goes to the upper end-point of $F$, is approximately a Generalised Pareto (GP) distribution $H_\gamma(\cdot\, / \sigma)$, which depends on a shape parameter $\gamma\in\Real$, called the extreme value index (EVI), and a scale parameter $\sigma>0$ \cite{balkema1974}.

The key benefits of this result are twofold: given a sample of independent and identically distributed (i.i.d.) observations from an unknown distribution $F$, the distribution $H_\gamma$ provides as an approximation for the distribution of rescaled excesses above a high threshold. By leveraging the quantile expression of $H_\gamma$, we derive an approximate formula for the extreme quantiles of $F$, a crucial tool in applications for assessing future risks.

The POT method is highly practical, yet its inferential theory remains complex (see \cite{dehaan+f06} for examples). Regardless of the inferential approach used, the theoretical analysis must account for the fact that the GP distribution is an inherently misspecified model for excesses, as the threshold must be fixed in practice. Additionally, the scale parameter is threshold-dependent, varying with the chosen threshold rather than being a fixed parameter as in classical statistical literature (e.g., \cite{vdv2000}). Moreover, the GP family is an irregular model, as the sign of the shape parameter influences its support, making likelihood-based inference—including Bayesian methods—notoriously challenging. In this paper, we develop inferential methods and the corresponding asymptotic theory while carefully addressing these complexities.

One of the most widely used statistical frameworks in applications is conditional inference, where the objective is to assess specific characteristics of a response variable $Y$ (such as its mean or quantiles) based on available information about covariates. In the second part of this work, we consider a response variable $Y$ following a generic distribution $F$ and a covariate vector $\bfX$ with law $P(\diff\bfx)$. Our focus is on modeling and analyzing the conditional distribution $F_{\bfx}$ of $Y$ given that $\bfX=\bfx$. Several methodologies have been proposed for modeling and statistically analyzing conditional extremes, including various conditional extreme value models and non- or semi-parametric regression approaches (see, e.g., \cite{resnick2014transition, goegebeur2014nonparametric}). In this work, we establish a link between the extremes of a response variable and its associated covariates by adopting a conditional distribution framework that complies with what we call the {\it proportional tail model}. Similar to the triangular array approach in \cite{einmahl2016}, this model allows the tail of $F_{\bfx}$ to vary according to a scale function, known as the {\it scedasis} function, while maintaining a constant extreme value index (EVI) (see also \cite{einmahl2022}). This approach enables the semi-parametric estimation of conditional tail probabilities and a non-parametric assessment of covariate effects using peaks above a high threshold and their associated {\it concomitant} covariates.

%
%
%
%
\subsection{Objectives and contributions}\label{sec:goals}
%
%
%
%
In the last decades, numerous inference methods for both unconditional and conditional extreme events have been developed (see, e.g., \cite{beirlant2006, daouia2022, dehaan+f06, einmahl2016, goegebeur2014, wang2015, wang2012}). While several studies have explored Bayesian inference for the POT method (see Ch. 11 in \cite{beirlant2006} and \cite{coles2003anticipating, do2016bayesian, fuquene2015semi, northrop2016, tancredi2006}), establishing a rigorous inferential theory remains a significant challenge. To the best of our knowledge, no asymptotic results have been derived for these Bayesian approaches.
More importantly, the ambitious and complex problem of predicting future extreme events within a proper statistical predictive framework has received little attention, with the notable exception of \cite{hall2002effect}. A practical and accessible approach to forecasting such events can be achieved through the Bayesian paradigm, which naturally yields the posterior predictive distribution—an estimator of the conditional distribution of an out-of-sample random variable, given that it exceeds a sufficiently high threshold, representative of future extreme events.
The main contributions of this article can be summarized as follows (with a detailed discussion to follow): (i) establishing a rigorous theoretical foundation for Bayesian inference in both unconditional and conditional settings within the POT framework; (ii) providing mathematical guarantees on the accuracy of forecasts based on the posterior predictive distribution.

In the first part of this paper, we develop the asymptotic theory for Bayesian inference within the classical POT framework. Specifically, we provide general and simple conditions on the prior distribution of the GP model's parameters under which we derive key results for the corresponding posterior distribution: consistency with a $\sqrt{k}$-contraction rate, the celebrated Bernstein-von Mises (BvM) theorem, and the asymptotic coverage probability of credible intervals. These results notably differ from those obtained in the block maximum context by \cite{padoan2022empirical}, as our approach allows for a broader selection of prior distributions and more general conditions. In particular, we accommodate families of informative proper priors for the shape parameter $\gamma$, and both non-informative improper priors and informative proper data-dependent priors for the scale parameter $\sigma$, as the prior specification for the latter is more nuanced. In contrast, \cite{padoan2022empirical} only considers the case of data-dependent priors. Moreover, while the theory in \cite{padoan2022empirical} relies on certain conditions on the density of $F$, which may be seen as restrictive, our results are derived under weaker, more standard conditions (see Ch. 2 in \cite{dehaan+f06}).

To develop a more comprehensive theory, we deeply investigate the frequentist properties of the empirical log-likelihood process, specifically in the context of the misspecified GP model. In this analysis, we derive its uniform convergence, the convergence of its derivatives (discussed in the supplement), as well as local and global bounds and a local asymptotically Gaussian expansion. Additionally, we obtain the contraction rates for the corresponding Maximum Likelihood Estimator (MLE).
Our results make several important contributions to statistical inference. First, they address a longstanding question: Is the MLE, computed over the entire parameter space, consistent, and is it the unique global maximizer of the likelihood? Second, our findings are essential for developing the asymptotic theory of Bayesian methods. While a version of the Bernstein-von Mises (BvM) theorem exists for misspecified models \citep{kleijn2012bernstein}, this result does not directly apply here. The GP distribution is an irregular statistical model that violates the smoothness conditions outlined in \citep{kleijn2012bernstein}, and the degree of misspecification in our setting is not fixed but instead varies with $n$. This makes the testability assumptions in \citep{kleijn2012bernstein} difficult to verify.

Beyond estimating the tail of $F$, practitioners are particularly interested in quantifying the quantiles of $F$ for exceptionally small exceedance probabilities, as these correspond to events more extreme than any observed so far. To address this, we extend the asymptotic theory of Bayesian methods by deriving consistency, contraction rates, asymptotic normality, and the asymptotic coverage probability of credible intervals for the posterior distribution of the so-called extreme quantiles (Ch. 3, 4 in \cite{dehaan+f06}). This is achieved through the development of a general Bayesian delta method, a widely applicable result that extends beyond EVT. The posterior distribution of extreme quantiles is a valuable tool for assessing the intensities of future extreme events, as it quantifies the uncertainty of their magnitude. However, to fully account for the uncertainty in predicting such events, we propose a genuine statistical predictive approach. We conclude the first part by introducing the posterior predictive distribution as an estimator for the conditional distribution of an out-of-sample random variable, given that it exceeds a sufficiently high threshold, representative of a future extreme event. We derive conditions under which this predictive distribution is Wasserstein consistent and quantify its contraction rate.

In the second part, we address the problem of Bayesian inference for the tail properties of the conditional distribution $F_{\bfx}$ of $Y$ given $\bfX=\bfx$. Building on the {\it tail proportionality} condition, we show that this inference can be achieved in two steps: first, Bayesian inference of the GP distribution parameters using peaks above a high threshold, and second, estimation of the scedasis function using the concomitant covariates. The first step is already described in the first part of the paper. For the second step, we specify a Dirichlet Process (DP) prior (e.g., Ch. 4.1 in \cite{ghosal2017}) for the conditional law of the concomitant covariate $\bfX$ given $Y > t$, which induces a prior on the scedasis function at $\bfx$. This function depends on the unknown marginal probability measure $P(\diff \bfx)$, which we treat as a nuisance parameter for simplicity, and we estimate it using two methods: a kernel-based method and a $K$-nearest neighbors approach. For the corresponding posterior distribution, we derive the same type of asymptotic theory established in the first part of the paper. With the kernel method, we also provide contraction rates for the posterior distribution of the scedasis function, which hold uniformly in $\bfx$. We then turn to the functional estimation of the marginal law of the concomitant covariates. Under the same Dirichlet Process prior, we show that the posterior distribution satisfies the BvM theorem over an infinite-dimensional Skohorod space. This result forms the basis for deriving a Kolmogorov-Smirnov-type statistical test to assess whether the concomitant covariates significantly affect the extremes of the response variable.

Extreme conditional quantiles are essential tools for assessing the risks associated with extreme events in a phenomenon, particularly when other concomitant dynamics reach certain levels. Under the proportional tail model, the posterior distribution of these quantiles is readily available, as it is induced by the posterior distributions of the GP distribution’s parameters and the scedasis function. To further refine our analysis, we also develop the asymptotic behavior of the posterior distribution. The final, but perhaps most significant, statistical problem addressed in this article is the forecasting within an extreme regression framework. In this context, we define the posterior predictive distribution as an estimator of the conditional distribution of an out-of-sample random variable, given that it exceeds an extreme conditional quantile and the concomitant covariates reach certain levels. For this predictive distribution, we establish minimal conditions to prove its Wasserstein consistency and quantify its contraction rate.

%
%
\subsection{Workflow}\label{sec:structure}
%
%
Section \ref{sec:pot_theory} presents key concepts and notation for the POT method, outlines the theoretical properties of the log-likelihood empirical process, develops the asymptotic theory for Bayesian inference within the POT framework, and establishes the Wasserstein contraction rates for the corresponding posterior predictive distribution. Section \ref{sec:ext_cond_Q} introduces the proportional tail model to link the extremes of a response variable to covariates, provides the posterior asymptotic theory for tail-related quantities in the conditional distribution, and derives the Wasserstein contraction rates for the associated posterior predictive distribution. Section \ref{sec:simulations} offers a comprehensive simulation study demonstrating the finite-sample performance of the proposed methodology. Section \ref{sec:application} concludes the paper with an application to real financial data, analyzing the change in crisis frequency over time. All proofs are provided in the supplement, which also includes additional results, details on posterior computations, simulations, and real data analysis.

%
%
\section{The Peaks-Over-Threshold approach}\label{sec:pot_theory}
%
%

%
%
\subsection{Background}\label{sec:background}
%
%

Consider a random variable $Y$ whose distribution $F$ is in the domain of attraction of a Generalised Extreme Value (GEV) distribution $G_\gamma$, in symbols $F\in\DoA(G_\gamma)$, where $\gamma\in\Real$ is the EVI that describes the  heaviness of distribution's tail \citep[][Theorem 1.1.3]{dehaan+f06}. This means that for any integer $m\geq1$, there are norming constants $a_m>0$ and $b_m\in \Real$ such that for all $y\in\Real$ that are continuity points of $G_\gamma$, $F^m(a_my+b_m)\to G_\gamma(y)$, as $m\to\infty$. Let $y^\star=\sup\{y:F(y)<1\}$ and $F_t$ be the {conditional distribution of $(Y-t)$ given that $Y>t$}. The domain of attraction condition (or first-order condition) can be equivalently formulated as follows. For $t<y^\star$, there is a scaling function $s(t)>0$ such that
\begin{equation}\label{eq:DoA}
	\lim_{t\uparrow y^{\star}} F_t(s(t)z)=H_\gamma(z),
\end{equation}
where $H_\gamma$ is a unit-scale GP distribution (\cite{balkema1974},  \cite[][Theorem 1.1.6]{dehaan+f06}). 
A possible choice for the norming constants is $b_m=U(m)$, where $U(t)=F^{\leftarrow}(1-1/t)$ with $F^{\leftarrow}(y)=\inf\{x:F(x)\geq y\}$ is the so-called {\it tail quantile}, $a_m=a(m)$ for a suitable positive function $a(\cdot)$ and for the scaling function one can set $s(t)=a(U^{\leftarrow}(t))$, \citep[][Ch. 1.2]{dehaan+f06}. In the sequel we consider this choice.
The GP is a family of two-parameters distributions defined as $H_{\bfvartheta}(z)=H_{\gamma}(z/\sigma)$ for all $z\in\Supp_{\bfvartheta}$, where $\bfvartheta=(\gamma,\sigma)^\top\in\Real\times(0,\infty)$,
$$
H_\gamma(z)=
\begin{cases}
	1-(1+\gamma z)_{+}^{-1/\gamma},& \text{if } \gamma\neq0,\\
	1-\exp(-z),& \text{if } \gamma=0,
\end{cases}
$$
with $(z)_+=\max(0,z)$, and $\Supp_{\bfvartheta}$ is $[0,\infty)$ if $\gamma\geq0$ while is $[0,-\sigma/\gamma]$ if $\gamma<0$. The GP density is $h_{\bfvartheta}(z)=h_{\gamma}(z/\sigma)/\sigma$, where
$$
h_{\gamma}(z)=
\begin{cases}
	(1+\gamma z)_{+}^{-(1/\gamma+1)},& \text{if } \gamma\neq0,\\
	\exp(-z), & \text{if } \gamma=0.
\end{cases}
$$
The log-likelihood of the density $h_\bfvartheta$, corresponding to a single observation, is defined for $z\geq0$ as
$$
\ell_\bfvartheta(z)=
\begin{cases}
	-\log \sigma -\left(1+\frac{1}{\gamma}\right)\log\left(1+\frac{\gamma}{\sigma}z\right),& \text{if } 1+\frac{\gamma}{\sigma}z>0,\\
	-\infty,& \text{otherwise}.
\end{cases}
$$
Observe that the log-likelihood is unbounded when $\gamma<-1$, more precisely for any $y\geq0$
$$
\lim_{\sigma\downarrow -\gamma z} \ell_\bfvartheta(z)=\infty.
$$
We recall that the Fisher information matrix corresponding to the GP log-likelihood is
\begin{equation}\label{eq:Info_matrix}
I_{\bfvartheta}=-\int_0^1 \frac{\partial^2 \ell_{\bfvartheta}}{\partial \bfvartheta \partial \bfvartheta^\top}
\left(
\frac{v^{-\gamma}-1}{\gamma}
\right)\diff v, 
\end{equation}
which is positive definite as soon as $\gamma>-1/2$. For this reason in the sequel we restrict the parameter space to $\Theta=(-1/2,\infty)\times(0,\infty)$.
%

%
%
\subsection{Empirical log-likelihood process asymptotics}\label{sec:preliminary}
%
%
%
In this section, we examine key frequentist properties of the empirical log-likelihood process associated with the GP density function. We derive its uniform convergence, along with that of its derivatives (detailed in the supplement), its local asymptotic expansion, and both local and global bounds. These results provide a comprehensive understanding of the GP likelihood theory, serving as a foundation for the asymptotic theory of Bayesian methods and offering valuable insights from a frequentist perspective. As a by-product, we establish the uniqueness and contraction rates of the MLE of the GP likelihood function for large samples, a new and noteworthy contribution.

For $n\geq1$, let  $(Y_1,\ldots,Y_n)$ be i.i.d. copies of $Y$ whose distribution satisfies $F_0\in\DoA(G_{\gamma_0})$.
The first-order condition \eqref{eq:DoA} implies that $\Prob(Y\leq y \mid Y>t)\approx H_{{\bfvartheta_0}}(y-t)$
for high enough threshold $t$, with $\bfvartheta_0=(\gamma_0,\sigma_0)^\top$
with $\sigma_0$ representative of $s_0(t)$ and $y\equiv y_t=t+s_0(t)z$, for all $z\geq0$.
Let $k=k(n)$ be {a} so-called {\it intermediate sequence}, with $k=o(n)$ and $k\to\infty$ as $n\to\infty$. 
A practical way of defining a high threshold is $t=U_0(n/k)$ and estimating $U_0(n/k)$ by the order statistic $Y_{n-k,n}$, where $Y_{1,n}\leq\cdots\leq Y_{n,n}$ are the $n$ order statistics, and $F$ by the empirical distribution $F_n$, we obtain $s_0(U_0(n/k))\approx a_0(F_n^{\leftarrow}(1-1/Y_{n-k,n}))= a_0(n/k)$. Accordingly, $\sigma_0$ is representative of $a_0(n/k)$. 
We focus on the normalised excess over-high threshold variables, also known as \lq\lq pseudo-observations'', 
\begin{equation*}\label{eq:pseudo_variables}
Z_i=\frac{Y_{n-i+1,n}-Y_{n-k,n}}{a_0(n/k)},\quad 1\leq i\leq k.
\end{equation*}
For an arbitrary Borel set $B$ the empirical probability measures relative to observed  unrescaled  peaks and pseudo-observations by $\Prob_n(B)=$ $k^{-1}\sum_{i=1}^k$ $\indic(Y_{n-k+i,n}-Y_{n-k,n} \in B)$ and $\mathbb{P}_n^{\text{pse}}(B)= k^{-1}\sum_{i=1}^k \indic(Z_i \in B)$, respectively.
Given a probability measure $P$ on a measurable space $(\mathcal{X}, \mathcal{B})$ and a measurable function $f:\mathcal{X}\mapsto \Real^p$ we use $Pf$ to denote $\int f\diff P$. 
Accordingly, $\Prob_nf=
k^{-1}\sum_{i=1}^k f(Y_{n-i+1,n}$ $-Y_{n-k,n})$ and $\Prob_n^{\text{pse}}f=
k^{-1}\sum_{i=1}^k f(Z_i)$.

Based on this, we define the empirical log-likelihood process relative to the observed unnormalized excesses over a high threshold as $\lik_n(\bfvartheta) = \Prob_n \ell_{\bfvartheta}$, $\bfvartheta \in \Theta$. Since $\sigma$ is not a fixed parameter and depends on the sample size $n$, to stabilize it as $n$ increases, we introduce the reparametrization $\bftheta := r(\bfvartheta) = (\gamma, \sigma/a_0(n/k))^\top$ for all $\bfvartheta \in \Theta$, yielding $\bftheta_0 = r(\bfvartheta_0) = (\gamma_0, 1)^\top$. Our theory is developed using the empirical log-likelihood process $L_n(\bftheta) = \Prob_n^{\text{pse}} \ell_{\bftheta}$, which is defined through the GP log-likelihood $\ell_{\bftheta}$. Note that $L_n(\bftheta) = \lik_n(\bfvartheta) + \log a_0(n/k)$. For convenience, we refer to $\lik_n$ and $L_n$ as the "realistic" and "theoretical" empirical log-likelihood processes, respectively. The former is the version used for practical inference, while the latter is employed to study the asymptotic properties of the process. Finally, we define the (theoretical) score and information processes as $S_n(\bftheta) = (\partial / \partial \bftheta) L_n(\bftheta)$ and $J_n(\bftheta) = (\partial^2 / \partial \bftheta \partial \bftheta^\top) L_n(\bftheta)$. In the following, given two vectors $\bfx, \bfy$ of equal size, $\bfx^\top \bfy$ denotes the usual vector product, while componentwise multiplication and division are denoted as $\bfx \bfy = (x_1 y_1, \dots, x_q y_q)^\top$ and $\bfx / \bfy = (x_1 / y_1, \dots, x_q / y_q)^\top$, respectively.

We now present some key results that are crucial for the asymptotic theory of the MLE and serve as the foundation for deriving the main findings in the subsequent section on the Bayesian approach. Let $B(\bftheta_0, \varepsilon) = {\bftheta \in \Theta : |\bftheta - \bftheta_0| < \varepsilon}$ denote the open ball centered at $\bftheta_0$ with radius $\varepsilon$, and let $B(\bftheta_0, \varepsilon)^\complement$ be its complement in $\Theta$. We define
$$
\widehat{\bftheta}_n{\in}\mathop{\mathrm{argmax}}_{\bftheta\in B(\bftheta_0,\varepsilon)} L_n(\bftheta)
$$
as a local maximizer of the empirical log-likelihood process. 
Note that the local maximiser of the realistic empirical log-likelihood process $\lik_n(\bfvartheta)$  satisfies $\widehat{\bfvartheta}_n=r^{-1}(\widehat{\bftheta}_n)$.
Note that a different  neighborhood of $\bftheta_0$ with compact closure in $\Theta$ might be considered for the definition of local MLE $\widehat{\bftheta}_n$, though $B(\bftheta_0,\varepsilon)$ is a natural choice. Note also that the first-order condition is equivalent to \citep[][Ch. 1]{dehaan+f06}
\begin{equation}\label{eq:first_order}
	\lim_{t\to\infty}\frac{U_0(ty)-U_{0}(t)}{a_0(t)} = \frac{y^{\gamma_0}-1}{\gamma_0},\quad \forall \, y>0.
\end{equation}
To control the asymptotic behavior of generic estimation procedures, 
we consider the following so-called second-order condition (e.g., Ch. 2 and Appendix B in \cite{dehaan+f06}). 
\begin{cond}\label{cond:second_order}
There is a rate (or second-order auxiliary) function $A$, i.e. a positive or negative function satisfying $A(t)\to 0$ as $t\to\infty$, such that:
\begin{inparaenum}
\item\label{cond:sec_order_formula} $A$ is regularly varying with index $\rho\leq 0$ (second order parameter) and
\begin{equation}\label{eq:second_order_par}
\lim_{t\to\infty}
\frac{\frac{U_0(tx)-U_{0}(t)}{a_0(t)} - \frac{x^{\gamma_0}-1}{\gamma_0}}
{A(t)}
=\int_1^x v^{\gamma_0-1}\int_1^v u^{\rho-1}\diff u \diff v.
\end{equation}
\item\label{cond_sec_order_bias} $\sqrt{k}A(n/k)\to\lambda\in\Real$ as $n\to\infty$.
\end{inparaenum}
\end{cond}
\begin{proposition}\label{prop:llik_local}
Under Condition \ref{cond:second_order} there exist $\varepsilon_0>0$ and constants $c_1,c_2,c_3>0$ such that the following three properties hold with probability tending to $1$ as $n\to\infty$: 
\begin{itemize}
\item When $\gamma_0>0$, 
$L_n$ is strictly concave on $B(\bftheta_0,\varepsilon_0)$ with a unique maximizer $\widehat{\bftheta}_n$ and 
\begin{equation}\label{eq:Taylor_positive}
-c_1\|\bftheta-\widehat{\bftheta}_n\|^2
\leq L_n(\bftheta) - L_n(\widehat{\bftheta}_n)\leq 
-c_2 \|\bftheta-\widehat{\bftheta}_n\|^2,\quad \bftheta\in B(\bftheta_0,\varepsilon_0).
\end{equation}
\item More generally, if $\gamma_0 > -1/2$, {$S_n(\bftheta_0) = O_{\Prob}(1/\sqrt{k})$ and}
\begin{equation}\label{eq:Taylor_negative}
L_n(\bftheta) - L_n(\bftheta_0) \leq (\bftheta-\bftheta_0)^\top {S_n(\bftheta_0)}-c_3 \|\bftheta-\bftheta_0\|^2,\quad \bftheta\in B(\bftheta_0,\varepsilon_0),
\end{equation}
\begin{equation}\label{eq:seconddev}
\sup_{\bftheta \in B(\bftheta_0, \epsilon_n)} \left\Vert
J_n(\bftheta)+ \bfI_{\bftheta_0} \right\Vert = o_{\mathbb{P}}(1), \text{ for any } \epsilon_n=R/\sqrt{k} \text{ with } R =o(\sqrt{k}/\log k \vee k^{1/2+\gamma_0^-}),
\end{equation}
where $\bfI_{\bftheta_0}$ is the Fisher information of the GP log-likelihood (defined as in \eqref{eq:Info_matrix} but with $\bftheta$ in the place of $\bfvartheta$) at $\bftheta_0$.

\end{itemize}
\end{proposition}
Next result uses Proposition~\ref{prop:llik_local} to establish the contraction rates
of the local MLE $\widehat{\bfvartheta}_n$.
\begin{corollary}\label{cor:consistency}
For any $R\to\infty$ satisfying $R=o(\sqrt{k})$ we have $\|\widehat{\bftheta}_n-\bftheta_0\|=O_{\mathbb{P}}(R/\sqrt{k})$. Accordingly, the normalised local MLE sequence $(\widehat{\bfvartheta}_n)_{n\geq 1}$ is $\sqrt{k}$-consistent, i.e.
\begin{equation}\label{eq:consistency}
\widehat{\gamma}_n=\gamma_0+O_{\bbP}(1/\sqrt{k})\quad\mbox{and}\quad \widehat{\sigma}_n=a_0(n/k)\big(1+O_{\bbP}(1/\sqrt{k})\big),
\end{equation}
and unique with probability tending to one, as $n\to\infty$.

\end{corollary}
\begin{rem}
Surprisingly, this result has not been established before. Previous works have mainly focused on the existence of a consistent local maximizer (potentially among several others) and the asymptotic behavior of the solutions to the likelihood equations, assuming they lie within a narrow neighborhood of the true parameter values. This neighborhood was defined by the conditions $|\gamma / (\sigma / a_0(n/k)) - \gamma_0| = O_{\mathbb{P}}(1/\sqrt{k})$ and $\sigma / a_0(n/k) = e^{O_{\mathbb{P}}(1)}$ \citep[e.g.,][Proposition 3.1]{d+f+d03}. Recently, \cite{einm22} showed that the log-likelihood is strictly concave within a slightly larger, but still shrinking, neighborhood of the form ${\bftheta : |\gamma - \gamma_0| + |\sigma / a_0(n/k) - 1| < R / \sqrt{k} }$, with probability tending to one as $n \to \infty$, ensuring that there is a unique maximizer within this region. However, this result does not rule out the possibility of other local maximizers outside this neighborhood. Corollary \ref{cor:consistency} addresses this gap, providing the contraction rates for a generic local maximizer. Nevertheless, it does not guarantee that the MLE, obtained by maximizing the likelihood over the entire parametric space, is consistent and coincides with the global unique likelihood maximizer.
\end{rem}
In Theorem \ref{thm1}, we derive global upper bounds for the theoretical empirical log-likelihood process, which serve as a crucial tool in answering the open question: Is the MLE computed over the entire parameter space consistent, and is it the unique global likelihood maximizer? The affirmative answer is given in Corollary \ref{cor:globalMLE}.
\begin{theorem}\label{thm1}
Under Condition \ref{cond:second_order} there exist $\varepsilon_0>0$ and constants $c_1,c_2>0$ and, for all large enough $\overline\tau>0$, there exist constants $c_3>0$, $c_4,c_5\in\mathbb{R}$ such that, with probability tending to $1$ as $n\to\infty$, 
\begin{align}
L_n(\bftheta) - L_n(\bftheta_0)&\leq (\bftheta-\bftheta_0)^\top S_n(\bftheta_0)-c_1\|\bftheta-\bftheta_0\|^2\quad \mbox{if }\bftheta\in B(\bftheta_0,\varepsilon_0), \label{eq:thm1-1}\\
L_n(\bftheta) - L_n(\bftheta_0) &\leq -c_2\quad  \mbox{if } \bftheta\in B(\bftheta_0,\varepsilon_0)^\complement,\label{eq:thm1-3}\\
L_n(\bftheta) - L_n(\bftheta_0) &\leq -\log\sigma-\frac{c(\tau)}{\sigma}-d(\tau),\quad \mbox{for all } \bftheta\in\Theta,\label{eq:thm1-4}
\end{align}
with $\tau=\gamma/\sigma$, $c(\tau)=c_3\mathds{1}_{\tau\leq \overline{\tau}}+(2\tau)^{-1}\log \tau \mathds{1}_{\tau> \overline{\tau}}>0$, $d(\tau)=c_4\mathds{1}_{\tau\leq \overline{\tau}}+(\log \tau+c_5)\mathds{1}_{\tau> \overline{\tau}}$.

\end{theorem}
\begin{corollary}\label{cor:globalMLE}
The MLE defined by $\widehat{\bftheta}_n\in\mathop{\mathrm{argmax}}_{\bftheta\in \Theta}L_n(\bftheta)$ is $\sqrt{k}$-consistent and, with probability tending to one, the maximiser is unique, i.e. $\widehat{\bftheta}_n$ is the unique global maximiser of the likelihood.
\end{corollary}
The global bounds presented in Theorem \ref{thm1} are essential for deriving contraction rates and the Bernstein-von Mises (BvM) theorem for the posterior distribution of the GP distribution's parameters (see the next section). In particular, the BvM result is derived from the asymptotic properties of the theoretical empirical log-likelihood process and its associated quantities (e.g., Ch. 7 in \cite{vdv2000}). We further extend the analysis of the theoretical empirical log-likelihood process by providing its local asymptotic expansion, from which we can deduce the asymptotic normality of the MLE. In the following, we denote a multivariate normal cumulative distribution by $\mathcal{N}(\bfmu, \bfSigma)$, where $\bfmu$ is the mean vector and $\bfSigma$ is the covariance matrix, which reduces to $\mathcal{N}(\mu, \sigma^2)$ in the univariate case. The corresponding probability measure is denoted as $\mathcal{N}(\cdot; \bfmu, \bfSigma)$.
\begin{proposition}\label{prop:consistency}
Assume that $F_0\in\DoA(G_{\gamma_0})$ 
and Condition \ref{cond:second_order} is satisfied. 
Then, for all fixed $c>0$ and $\epsilon_n=c/\sqrt{k}$ we have 
\begin{equation*}
\sup_{\bftheta \in B(\bftheta_0,  \epsilon_n)} 
\left|
L_n(\bftheta)-L_n(\bftheta_0)-(\bftheta-\bftheta_0)^\top S_n(\bftheta_0)+\frac{1}{2k}(\bftheta-\bftheta_0)^\top \bfI_{\bftheta_0} (\bftheta-\bftheta_0)
\right|=o_{\Prob}{\left(
\frac{1}{k}
\right)}.
\end{equation*}
In particular, $\sqrt{k}S_n(\bftheta_0)\stackrel{d}{\rightarrow}\mathcal{N}(\lambda {\boldsymbol{\mu}}, {\boldsymbol{V}})$, where $\lambda{\boldsymbol{\mu}}$ is a bias term { and }
\begin{equation}
{	
\boldsymbol{V}=
\left(
\begin{matrix}
	\frac{5 \gamma_0^2 + 6\gamma_0+2}{(1+2\gamma_0)^2(1+\gamma_0)^2}, 
	& 
	\frac{1+\gamma_0}{(1+2\gamma_0)^2}
	\\
	\frac{1+\gamma_0}{(1+2\gamma_0)^2},
	& \frac{(1+\gamma_0)^2}{(1+2\gamma_0)^2}
\end{matrix}
\right),
}
\end{equation}
see Section {2.7.1} of 
supplement for details.
Accordingly,  {with $\boldsymbol{b}=\bfI_{\bftheta_0}^{-1}\boldsymbol{\mu}$ and $\boldsymbol{\Sigma}= \bfI_{\bftheta_0}^{-1} \boldsymbol{V} \bfI_{\bftheta_0}^{-1}$,}
$$
\sqrt{k}(\widehat{\bftheta}_n-\bftheta_0)\stackrel{d}{\rightarrow}\mathcal{N}(\lambda  \bfb,{\boldsymbol{\Sigma}} ).
$$
\end{proposition}
\begin{rem}\label{rem:spatial_case}
Recently, \cite{einm22} derived similar asymptotic results for the GP log-likelihood in a more complex nonstationary space-time framework, assuming no bias. In contrast, we work in a simpler setup but derive the theory under the more general assumption that bias may exist. While the asymptotic bias and variance of the MLE in our study match those reported in \cite[][Theorem 3.4.2]{dehaan+f06}, our result is the first to establish asymptotic normality for the global likelihood maximiser.
\end{rem}
%

%
%
\subsection{Asymptotic theory of the posterior distribution}\label{sec:posterior_POT}
%
%

In this section we study the asymptotic properties of a Bayesian procedure for inference with the POT approach. 
We assume to work with a prior distribution on $\bfvartheta \in \Theta$ {with density} of the following form
\begin{equation}\label{eq:prior_ist}
\pi(\bfvartheta)= \pi_{\text{sh}}(\gamma) \pi_{\text{sc}}^{(n)}
\left( \sigma\right), \quad \bfvartheta \in \Theta,
\end{equation}
where $ \pi_{\text{sh}}$ is a prior density on $\gamma$ and for each $n=1,2,\ldots$, $\pi_{\text{sc}}^{(n)}$ is a prior density on $\sigma$, whose expression may or may not depend on $n$.
Although the prior in \eqref{eq:prior_ist} assumes independence between $\gamma$ and $\sigma$, it enables the specification of fairly flexible forms of their joint density. To establish our main results on the posterior distribution of these parameters, we need to work with a (genuine or empirical) prior density that satisfies the following weak conditions.
\begin{cond}\label{cond:prior}
The densities  $\pi_{\text{sh}}$ and $\pi_{\text{sc}}^{(n)}$ are such that:
\begin{inparaenum}
	\item\label{cond:pisc} For each $n=1,2,\ldots$, $\pi_{\text{sc}}^{(n)}:\Real_+\to\Real_+$ and
	\begin{inparaenum}
		\item[(a.1)]\label{posit} there is a constant $\delta>0$ such that $\pi_{\text{sc}}^{(n)}(a_0(n/k))a_0(n/k)>\delta$ and for any constant $\eta>0$ there is $\epsilon>0$ such that 
$$
\sup_{ 1-\epsilon < \sigma <1+\epsilon} \left|
\frac{\pi_{\text{sc}}^{(n)}(a_0(n/k)\sigma)}{\pi_{\text{sc}}^{(n)}(a_0(n/k))}-1 \right|	<\eta;
$$
\item[(a.2)]\label{piscbound} there is $C>0$ such that  ${\sup_{\sigma >0} \sigma a_0(n/k) \pi_{\text{sc}}^{(n)}(a_0(n/k)\sigma) \leq C}$;
\end{inparaenum}
Inequalities (a.1)-(a.2) hold with probability tending to $1$, for fixed $\delta, \eta, \epsilon, C$, if $\pi_{\text{sc}}^{(n)}$ is data-dependent.
%
\item\label{cond:pish}  $\pi_{\text{sh}}$ is a positive and continuous function at $\gamma_0$ such that: $\int_{-1/2}^0 \pi_{\text{sh}}(\gamma)\,\mathrm{d}\gamma<\infty$, $\sup_{\gamma>0}\pi_{\text{sh}}(\gamma)<\infty$.
\end{inparaenum}
\end{cond}
Below we provide concrete examples of informative and improper non-informative prior distributions that satisfy such conditions.
\begin{example}[Informative data dependent prior]\label{ex:data_dependent}
Let $\pi_{\text{sc}}^{(n)}(\cdot)=\pi(\cdot/\widehat{\sigma}_n)/\widehat{\sigma}_n$, where $\pi$ is an informative prior density on $(0,\infty)$ and $\widehat{\sigma}_n$ is an estimator of $a_0(n/k)$. Then,
$$
\pi_{\text{sc}}^{(n)}(\sigma a_0(n/k))a_0(n/k)=\pi(\sigma a_0(n/k)/\widehat{\sigma}_n) a_0(n/k)/\widehat{\sigma}_n.
$$
Now, if $a_0(n/k)/\widehat{\sigma}_n\stackrel{\bbP}{\rightarrow}1$ as $n\to\infty$ (examples are the MLE, generalised probability weighted moment estimator, etc.) and $\pi$ is continuous and $\sigma\mapsto\pi(\sigma t)$ is uniformly integrable  { in $\sigma$ for} $t$ in a neighbourhood of $1$ (examples are {G}amma, {I}nverse-{G}amma, Weibull, Pareto, etc.), then Condition \ref{cond:prior}\ref{cond:pisc} is satisfied. The informative joint prior density is completed setting $\pi_{\text{sh}}(\gamma)=\pi(\gamma)$, where $\pi$ is a  probability density on ${(}-1/2,\infty)$ assumed continuous and bounded away from infinity. 
\end{example}
\begin{example}[Non-informative improper prior]\label{ex:uninformative}
Consider a uniform distribution on $\log\sigma$ so that $\pi_{\text{sc}}^{(n)}(\sigma)=\pi(\sigma)\propto 1/\sigma$, for any given $n\geq 1$. Accordingly, well known non-informative prior densities are: the uniform prior $\pi(\bfvartheta)\propto \sigma^{-1}$, the maximal data information $\pi(\bfvartheta)\propto \sigma^{-1}\exp{-(\gamma+1)}$ and the Jeffreys prior $\pi(\bfvartheta)\propto \sigma^{-1}((1+\gamma)(1+2\gamma)^{1/2})^{-1}$, with $\sigma>0$ and $\gamma>-1/2$ (e.g., \cite{northrop2016}). In these cases
$$
\pi(a_0(n/k))a_0(n/k)=1,\quad \frac{\pi(a_0(n/k)\sigma)}{\pi(a_0(n/k))}=\frac{1}{\sigma},
$$
and so Condition \ref{cond:prior}\ref{cond:pisc} is trivially satisfied.
\end{example}
Given a prior density $\pi$ (proper informative or improper non-informative), accordingly, the posterior distribution on the parameters $\bfvartheta$ of the GP distribution is defined as
\begin{equation}\label{eq:posterior_pot}
\Pi_n(B)=\frac{\int_{B}\exp(k\lik_n(\bfvartheta))\pi(\bfvartheta)\diff \bfvartheta}
{\int_{\Theta}\exp(k\lik_n(\bfvartheta))\pi(\bfvartheta)\diff \bfvartheta},
\end{equation}
for all measurable sets $B\subset \Theta$.
Consistently with the frequentist context, the asymptotic  theory of the posterior distribution is derived working with the theoretical empirical log-likelihood process $L_n$. 
The corresponding posterior distribution is denoted by ${\Upsilon_n=\Pi_n}\circ r^{-1}$ (see Section \ref{sec:preliminary}).
 \begin{theorem}\label{theo:main_posterior}
Assume $F_0\in\DoA(G_{\gamma_0})$  
and that Condition \ref{cond:second_order} is satisfied and that the prior density {$\pi(\boldsymbol{\vartheta})$, $\boldsymbol{\vartheta}\in\Theta$}, satisfies Condition \ref{cond:prior}. Then: 
\begin{itemize}
\item (Contraction rate) {$\Upsilon_n$} is consistent with $\sqrt{k}$-contraction rate, namely there is a $R>0$ such that for all sequences $\varepsilon_n\to0$, satisfying  $\sqrt{k}\varepsilon_n\to \infty$ as $n\to\infty$,
$$
{\Upsilon_n}\left(B(\bftheta_0, \varepsilon_n)^\complement\right)= O_\bbP\left(\exp\left(-Rk\varepsilon_n^2\right)\right).
$$
\item (Bernstein-von Mises) As $n\to\infty$, 
$$
\sup_{B\subset \Theta}|\Upsilon_n(\{\bftheta:\sqrt{k}(\bftheta-\bftheta_0)\in B\})-\mathcal{N}(B;\bfI_{\bftheta_0}^{-1}\sqrt{k}S_{n}(\bftheta_0), \bfI_{\bftheta_0}^{-1})|=o_\bbP(1),
$$
where $B$ represents any Borel set in $\Theta$.
\item {(Coverage probability)} If $\lambda=0$, for any  $\alpha\in(0,1)$, as $n\to\infty$,
$$
\bbP\left(
{\left\{
\gamma_0 \in  (\Upsilon_{n;\text{sh}}^{\leftarrow}(\alpha/2), \Upsilon_{n;\text{sh}}^{\leftarrow}(1-\alpha/2))
\right\}
}	
\right)=1-\alpha+o(1), \quad 
$$
where $\Upsilon_{n;\text{sh}}^{\leftarrow}(1-a)$ is the $(1-a)$-quantile of the { posterior distribution of $\gamma$.}
\end{itemize}
\end{theorem}
Note that Theorem \ref{theo:main_posterior} implies that the posterior distribution $\Pi_n$ and its marginal distributions asymptotically concentrate around the true parameters, and they shape as a normal distribution. Moreover, the coverage probabilities of credible intervals for $\gamma_0$ achieve asymptotic  the nominal level. Beyond estimating $\gamma_0$ to assess the heaviness of the data distribution’s tail, the primary objective of EVT is to infer extreme events beyond those observed in the past. This can be achieved by examining the $(1-p)$-quantile of $F_0$ for an exceedance probability $p = p_n$ such that $np \to c \geq 0$ as $n \to \infty$ (or, more generally, $p = o(k/n)$). The key case occurs when $c \leq 1$, implying that in a sample of size $n$, at most one observation is expected to exceed such an extreme quantile. Since $F_0$ is unknown in applications, we can use the approximation $1-p = F_0(y_p) \approx 1 - (k/n) (1 - H_{\boldsymbol{\vartheta}_0}(y_p - U_0(n/k)))$ for large $n$ (derived from the first-order condition \eqref{eq:DoA}), and $y_p = F_0^{-1}(1-p)$, to obtain as $n \to \infty$ the quantile approximation
%
%
\begin{equation}\label{eq:extreme_quantile}
	F_0^{\leftarrow}(1-p)\approx U_0\left(\frac n k \right) + H^{\leftarrow}_{\boldsymbol{\vartheta}_0}\left(1-\frac{np}{k}\right)= U_0\left(\frac n k \right) + \sigma_0\frac{\left(\frac{k}{np}\right)^{\gamma_0}-1}{\gamma_0},
\end{equation}
\citep[e.g., ][Ch. 3]{dehaan+f06}. For each given $n$, $k$ and $p$, the right hand-side of \eqref{eq:extreme_quantile} is a continuous map $T_n:\Theta\to\Real$  and therefore $\Pi_n$ induces a posterior distribution $\Pi_n\circ T_n^{-1}$ on the approximate extreme quantile. When studying its properties we take into account that $U_0(n/k)$ is not covered by the Bayesian procedure but it is actually frequentistically estimated by $Y_{n-k,n}$.  We 
refer to the map, the posterior distribution and the extreme quantile as $\widetilde{T}_n$, $\widetilde{\Pi}_n:=\Pi_n\circ \widetilde{T}_n^{-1}$ and $Q(p)=Y_{n-k,n}+H^{\leftarrow}_{\boldsymbol{\vartheta}}(1-np/k)$, respectively.

We introduce a Bayesian delta method, which is useful for deriving the asymptotic theory of the posterior distribution of extreme quantiles. To the best of our knowledge, it has not been previously explored and is of independent interest, as it ca be applied beyond the EVT context in other settings.
%
\begin{theorem}\label{theo:Bayesian_delta_method}
{Consider a statistical model with parameter}
${\bftheta} \in\Theta\subset \Real^m$. 
Let $\pi(\bftheta)$ be a prior density on $\bftheta$, which can possibly be data dependent as described e.g. in \eqref{eq:prior_ist}, and $\Pi_n(\bftheta\in\cdot)$ be the posterior distribution of $\bftheta$. Assume the following conditions:
\begin{inparaenum}
\item\label{cond:posterior} Let $\bfv_n{\in(0,\infty)^m},\bftheta_n\in\Real^m$ be sequences such that $\bfv_n\to{\boldsymbol{\infty}}$ and $\bftheta_n\stackrel{\bbP}{\rightarrow} \bftheta_0$ and 
$$
\sup_{B\subset\Theta}|\Pi_n(\{\bftheta:\bfv_n(\bftheta-\bftheta_n)\in B\})-\mathcal{N}(B; {\boldsymbol{0}}, {\boldsymbol{D}})|=o_{\bbP}(1),
$$
where $B$ represents any measurable subset of $\Theta$ {and $\boldsymbol{D}$ is positive definite}.
\item\label{cond:transformation}  
Let $l \leq m$ and  $T_n:\Theta \to \Real^l$ be a sequence of continuously differentiable maps such that, 	for a sequence $\bfw_n\in \mathbb{R}^l$ and a $(l\times m)$ matrix $\bfJ$ of full rank, we have $\nabla \overline{T}_n(\bfDelta)\to\bfJ$ uniformly on compact sets as $n\to\infty$, where
$$
\overline{T}_n(\bfDelta)=\bfw_n(T_n(\bftheta_n+\bfv_n^{-1}\bfDelta)-T_n(\bftheta_n)), \quad \bfDelta\in\Real^l.
$$
\end{inparaenum}
Then, 
$$
\sup_{B\subset {\Real^l}}|\Pi_n(\{\bftheta: \bfw_n(T_n(\bftheta)-T_n(\bftheta_n))\in B\})-\mathcal{N}(B;0,\bfJ \boldsymbol{D} \bfJ^\top)|=o_{\bbP}(1),
$$
where $B$ represents any measurable subset of $ \Real^l$.
\end{theorem}

Finally, the following corollary presents the asymptotic properties of the posterior distribution for extreme quantiles, derived directly from Theorems \ref{theo:main_posterior} and \ref{theo:Bayesian_delta_method}.
In view of the next result, we introduce the following function (e.g., Ch. 4.3 in \cite{dehaan+f06}) for any $t>1$ and $\gamma{>-1/2}$,
$$
q_{\gamma}(t):=\frac{\partial H_{\gamma}^{\leftarrow}}{\partial \gamma}(1-1/t)=\int_{1}^t v^{\gamma-1} \log v \diff v.
$$
\begin{corollary}\label{theo:quantile_posterior}
Assume that the conditions of Theorem \ref{theo:main_posterior} are satisfied and, in the special subcase where $\rho=0$, further  assume that $\gamma_0<0$. 
Then, for $p=o(k/n)$ such that $\log(k/np)=o(\sqrt{k})$:
\begin{itemize}
	\item (Contraction rate) There is a $R>0$ such that for all sequences $\varepsilon_n\to0$, satisfying  $\sqrt{k}\varepsilon_n\to \infty$ and $\varepsilon_n\log(k/np)\to 0$ as $n\to\infty$,
	$$
	\widetilde{\Pi}_n\left(\left\{Q(p)\in \Real: \left|\frac{Q(p)-F_0^{\leftarrow}(1-p)}{q_{\gamma_0}(k/(np))a_0(n/k)}\right|>\varepsilon_n\right\}\right)= O_\bbP\left(\exp\left(-Rk\varepsilon_n^2\right)\right).
	$$
	\item (Bernstein-von-Mises) For $v_n=\sqrt{k}/(q_{\gamma_0}(k/(np))a_0(n/k))$, as $n\to\infty$
	$$
	\sup_{B\subset \Real}\left|\widetilde{\Pi}_n\left(\left\{Q(p)\in \Real:v_n\left(Q(p)-F_0^{\leftarrow}(1-p)\right)\in B\right\}\right)-\mathcal{N}(B;\widetilde{\Delta}_n, \widetilde{V})\right|=o_\bbP(1),
	$$
	where $B$ is any measurable subset of $\Real$, see Section 3.7 and Equation (3.32) of the supplement for the explicit expression of $\widetilde{\Delta}_n$ and $\widetilde{V}$.
	\item (Coverage probability) {If $\lambda=0$}, for any fixed $\alpha\in(0,1)$, as $n\to\infty$
	$$
	\bbP\left(\left\{F_0^{\leftarrow}(1-p)\in (\widetilde{\Pi}_n^{\leftarrow}(\alpha/2), \widetilde{\Pi}_n^{\leftarrow}(1-\alpha/2))\right\}\right)= 1-\alpha+o(1)
	$$
	%
where $\widetilde{\Pi}_n^{\leftarrow}(1-a)$ is the $(1-a)$-quantile of $\widetilde{\Pi}_n$, for $a\in(0,1)$.
\end{itemize}
\end{corollary}

\begin{rem}
The restriction $\gamma_0<0$, imposed in the special case where Condition \ref{cond:second_order} holds with $\rho=0$, helps to control bias when applying formula \eqref{eq:extreme_quantile}. This assumption is also common in the frequentist framework; \citep[e.g.][Theorem 4.3.1]{dehaan+f06}.
%
\end{rem}

%
%
\subsection{Asymptotic theory of predictive distribution}\label{sec:predictive_POT}
%
%

The posterior distribution of extreme quantiles is undoubtedly a valuable tool for addressing the challenging task of assessing yet-to-occur extreme events, as it inherently provides a measure of uncertainty. However, at its core, it remains a method for inferring a distributional parameter. A more comprehensive approach to quantifying uncertainty in forecasting future extreme events involves employing a genuine statistical predictive framework.

A practical way to achieve statistical prediction is through the posterior predictive distribution. Given a past sample $\bfY_n=(Y_1,\ldots,Y_n)$, we consider an independent out-of-sample random variable $Y^{\star}$ as a representative of a future event. We then focus on the conditional distribution
$F_{0,n}^{\star}(y)=\bbP(Y^{\star}\leq y \mid Y^{\star}>U_0(n/k),\bfY_n)$,
which we refer to as the predictive distribution of an extreme event. While conditioning on $\bfY_n$ is technically unnecessary—since $Y^{\star}$ and $\bfY_n$ are independent—we retain it to emphasize the crucial role of past data in defining its estimator.
Following the Bayesian inferential framework outlined in Section \ref{sec:posterior_POT}, a natural estimator of the predictive distribution of an extreme event is given by the posterior predictive distribution.
\begin{equation}\label{eq:predictive_dist}
\widehat{F}_n^{\star}(y)=\int_{\Theta} H_{\boldsymbol{\vartheta}}(y-Y_{n-k,n})
\Pi_n(\diff \bftheta).
\end{equation}
The posterior predictive distribution can serve as a powerful tool for forecasting future extreme peaks. For instance, one can obtain a point forecast by computing the quantile $\widehat{F}_n^{{\star}\leftarrow}(1-p^{\star})$ for a small $p^{\star}\in(0,1)$, or derive a more comprehensive prediction by identifying an entire region of plausible future values based on the highest posterior predictive density (e.g., \cite{robert2007}).
Given the significant societal impact of extreme events, it is crucial to assess the accuracy of the proposed forecasting method. The next result establishes that  forecasts derived from our posterior predictive distribution are asymptotically reliable, as the latter approaches to the true predictive distribution as the sample size grows.
To quantify the closeness between two distributions $F$ and $G$, we use the Wasserstein distance of order $v$ for $v\geq 1$, i.e. $W_v(F,G)=(\int_0^1 |F^{\leftarrow}(p)-G^{\leftarrow}(p)|^v\diff p)^{1/v}$.
Moreover, we recall that by the scaling property of the Wasserstein distance, for $a_0(n/k)>0$, we have
$$
W_v(\widehat{F}_n^{\star}, F_{0,n}^{{\star}})=a_0(n/k)W_v(\widehat{F}_n^{\star}(a_0(n/k)\cdot),F_{0,n}^{\star}(a_0(n/k)\cdot)).
$$
\begin{theorem}\label{theo:pred_wass}
Assume that the conditions of Theorem \ref{theo:main_posterior} are satisfied. Assume also that Condition \ref{cond:prior}\ref{cond:pish} changes as: $\pi_{\text{sh}}$ is positive and continuous at $\gamma_0$ and there is $B\subseteq(-1/2,1/v)$ such that $\pi_{\text{sh}}(\gamma)=0$, for all $\gamma\notin B$, and 
$$
\int_{B \cap (-1/2,0)} \pi_{\text{sh}}(\gamma)\diff \gamma <\infty, \quad
\int_B (1-\gamma v)^{-1/v}\pi_{\text{sh}}(\gamma)\diff\gamma <\infty.
$$
Then, for all sequences $\varepsilon_n\to0$, satisfying  $k\varepsilon_n^2/\log(k)\to \infty$ as $n\to\infty$, we have
$$
\frac{W_v(\widehat{F}_n^{\star}, F_{0,n}^{\star})}{a_0(n/k)}=O_\bbP(\varepsilon_n).
$$
\end{theorem}
\begin{rem}\label{rem:extra_cond}
The alternative condition included in Theorem \ref{theo:pred_wass} is a requirement to work with the Wasserstein metric of order $v$, Specifically, ensuring the integrability of the $v$-moment of the GP distribution with respect to the prior density $\pi_{\text{sh}}$ is necessary for the theory to hold. This condition is relatively weak and is naturally satisfied by several common prior distributions. Examples include: Uniform prior on $(-1/2,1/v-\varepsilon)$ for an arbitrary small $\varepsilon\geq 0$; A Beta prior on the transformed parameter $(1-\gamma v)/(1+v/2)$, where the Beta shape parameters are $(\alpha+1/v, \beta)$, for $\alpha,\beta>0$.
\end{rem}
\section{Extreme regression}\label{sec:ext_cond_Q}
In this section, we introduce a Bayesian framework for inference and prediction of extremes of a response variable that is linked to some covariates that we model through the {\it proportional tail model} for conditional extremes. Similar to the triangular array approach in \cite{einmahl2016}, our key assumption is that the upper tail of the conditional distribution of the response, given the covariates, changes according to a scaling factor, while the EVI remains unchanged.
Our method leverages peaks above a high threshold, along with the corresponding concomitant covariates, to estimate the marginal tail probability parametrically and the effect of covariates non-parametrically. This foundation enables semi-parametric inference on extreme conditional quantiles and facilitates forecasting of conditional extremes.
To the best of our knowledge, this is the first approach that jointly models and estimates both the marginal extremal properties of the response and the dependence structure induced by covariates. 

\subsection{Proportional tail model}\label{sec:proportional_tail}
Let $(Y, \bfX)$ be a random vector on $\Real\times[0,1]^d$. We denote the marginal distribution of $Y$ by $F_0$, the marginal law of $\bfX$ by $\Pmu_0$, the conditional distribution of {$Y$ given that $\bfX=\bfx$} by $F^{(0)}_{\bfx}$ and the corresponding $(1-p)$-quantile by $F^{(0)\leftarrow}_{\bfx}(1-p)$. 
We assume that $F_0$ is absolutely continuous and that satisfies $F_0\in\DoA(G_{\gamma_0})$ and  the conditional distribution $F^{(0)}_{\bfx}$ satisfies the {\it tail proportionality} condition, 
i.e. there is a positive {bounded} function $c_0$ on $[0,1]^d$, named {\it scedasis function} \citep{einmahl2016},  such that
\begin{equation}\label{eq:tail_prop}
	\lim_{y\to y^*} \frac{1-F^{(0)}_{\bfx}(y)}{1-F_0(y)}=c_0(\bfx), \quad \bfx\in[0,1]^d.
\end{equation}
Note that Condition \eqref{eq:tail_prop} together with the first-order condition in \eqref{eq:first_order} implies
$$
\lim_{t\to\infty}\frac{U^{(0)}_{\bfx}(ty)-U^{(0)}_{\bfx}(t)}{a_0(t)} = (c_0(\bfx))^{\gamma_0}\frac{y^{\gamma_0}-1}{\gamma_0},\quad \forall \, y>0,
$$
{where $U_{\boldsymbol{x}}^{(0)}(t)=F^{(0)\leftarrow}_{\bfx}(1-1/t)$}.
This means that 
the conditional distributions $F_{\bfx}^{(0)}$ changes 
according to the scaling function $(c_0(\bfx))^{\gamma_0}$, while the heaviness of its tail remains unchanged, since its tail index is equal to $\gamma_0$ no matter what is $\bfx$.

Under this framework {and assuming convergence in  \eqref{eq:tail_prop} to be uniform,} we obtain the following asymptotic approximations.  First,  according to the univariate case, 
for all $z>0$ and with  $y=U_0(n/k)+{z}$ we have $\Prob(Y\leq y \mid Y>U_0(n/k))\approx H_{\bfvartheta_0}(y-U_0(n/k))$ as $n\to\infty$.
Leveraging on this and on \eqref{eq:tail_prop} we obtain that for all measurable $B\subset [0,1]^d$ 
\begin{eqnarray}
	\nonumber \bbP(Y>y, \bfX\in B)&=&\int_B (1-F_{\bfx}^{(0)}(y))\Pmu_0(\diff\bfx)\\
	\nonumber &{\approx}&\int_B c_0(\bfx)(1-F_0(y))\Pmu_0(\diff\bfx)\\
	\label{eq:joint_tail}&{\approx}&\frac{k}{n}\Pmu^{*}_0(B)(1-H_{\bfvartheta_0}(y-U_0(n/k))),
\end{eqnarray}
where $\Pmu^{*}_0(\diff\bfx):=c_0(\bfx)\Pmu_0(\diff\bfx)$ and the approximations in the last two lines hold for $n\to\infty$.  
The above result entails that the conditional distribution of $\bfX$ given that $Y>U_0(n/k)$ is asymptotically approximated by the probability measure $\Pmu_0^*$, as $n\to\infty$ (see Lemma 3.5 in the supplement).
Second, as a direct consequence, the conditional tail  probability $\bbP(Y>y \mid X\in B)$ can be approximated in turn by the right-hand side of \eqref{eq:joint_tail} divided by $\Pmu_0(B)$, as $n\to\infty$. As a result one obtains 
for the conditional distribution $\bbP(Y\leq y \mid X\in B)=1-\bbP(Y>y \mid X\in B)$ a useful approximation that we refer to as the conditional proportional tail model. Third, the result
\begin{eqnarray}
	\nonumber & \bbP(&Y\leq y, \bfX\in B \mid Y>U_0(n/k))\\
	\label{eq:joint_dist_yx} &\approx& \Pmu_0^*(B) - \Pmu_0^*(B)(1-H_{\bfvartheta_0}(y-U_0(n/k)))
	=\Pmu_0^*(B)H_{\bfvartheta_0}(y-U_0(n/k)),
\end{eqnarray}
suggests that the joint conditional distribution of $(Y, \bfX)$ given that $Y>U_0(n/k)$ factorises asymptotically to a product of marginal distributions, which is useful in the next section to derive a Bayesian procedure for the inference on the parameters $(\bfvartheta_0, \Pmu^*_0)$, that allows in turn to make inference about $\bbP(Y\leq y \mid X\in B)$, through 
 the conditional proportional tail model. In the next section we study functional estimation of $\Pmu^*_0(B)$, namely for an infinite collection of sets $B$, as it allows to perform hypothesis testing to verify the effect of the covariates $\bfX$ on $Y$, given that $Y>U_0(n/k)$,
and both pointwise and functional estimation of its density $\diff \Pmu^*_0(\bfx)/\diff \Pmu_0(\bfx)$.
Note that $\Pmu^{*}_0(B)/\Pmu_0(B)\approx c_0(\bfx)$ when $B \downarrow\{\bfx\}$, and the conditional proportional tail model approximation of conditional distribution becomes
\begin{equation}\label{eq:tail_probability}
	\bbP(Y\leq y| \bfX=\bfx)\approx 1- c_0(\bfx)\frac{k}{n}\left(1+\gamma_0 \frac{y-U_0(n/k)}{a_0(n/k)}\right)_+^{-1/\gamma_0},
\end{equation}
as $n\to\infty$. In this context, the aim is to infer the scedasis function $c_0$ and more importantly for applications, the extreme quantiles of the conditional distribution.   Assuming that {$p=o(k/n)$} as $n\to\infty$, then exploiting the right-hand side of formula \eqref{eq:tail_probability} one obtains, for any $\bfx \in[0,1]^d$, the following approximation for the $(1-p)$-quantile $F^{(0) \leftarrow}_\bfx(1-p)$ of the conditional distribution,
\begin{equation}\label{eq:conditional_quantile}
	F^{(0)\leftarrow}_\bfx(1-p) \approx U_0(n/k) + H^{\leftarrow}_{\bfvartheta_0}\left(1-\frac{np}{k}\frac{1}{c_0(\bfx)}\right),
\end{equation}
as $n\to \infty$. Therefore, inference on $F^{(0)\leftarrow}_\bfx(1-p)$ can be achieved 
{leveraging that }on $(c_0, \bfvartheta_0)$.
%
%
%
\subsection{Asymptotic theory of the posterior distribution}\label{sec:posterior_EQ}
%
%
%
Let $(Y_i,\bfX_i)_{1\leq i \leq n}$ be a sample of i.i.d. copies of $(Y, \bfX)$. The Bayesian inference for the proportional tail model and related quantities is grounded on the joint statistical model $\{\mathcal{H}_{\bftheta}^k \times {\Pmu^*}^k, \bfvartheta\in\Theta, \Pmu^*\in \Psp\}$, which is motivated by the approximation \eqref{eq:joint_dist_yx}. In particular, $\mathcal{H}_{\bftheta}^k$ and  ${\Pmu^*}^k$ are the probability measures of  $k$ independent GP variables and concomitant covariates. Moreover, $\Psp$ is the family of Borel probability measures on $[0,1]^d$. The model is fitted to
the subsample $(Y_{n-i+1,n}-Y_{n-k,n},\bfX_{n-i+1,n})_{1\leq i\leq k}$ of peaks $(Y_{n-i+1,n}-Y_{n-k,n})_{1\leq i\leq k}$ above a high threshold $Y_{n-k,n}$ and concomitant covariates $(\bfX_{n-i+1,n})_{1\leq i\leq k}$, 
with $\bfX_{1,n},\ldots, \bfX_{n,n}$  that are the covariates associated to the order statistics {$Y_{1,n}<\cdots< Y_{n,n}$}. 
Note that the continuity of the distribution $F_0$ ensures that there are almost surely no ties. 
In this way there are several sources of misspecification: the exceedances  are dependent and only approximately distributed according to the Pareto distribution $H_{\bfvartheta_0}$, the corresponding concomitant covariates are only approximately distributed according to the law $\Pmu^*_0$, exceedances  and concomitant covariates are dependent and only   approximately independent from each other.  Despite that, we can show that the posterior distribution of the proportional tail model parameters and the conditional extreme quantiles enjoy good asymptotic properties.

We specify the prior distribution for the proportional tail model parameters as ${\Lambda}(\diff\bfvartheta, \diff \Pmu^*)=\Pi(\diff\bfvartheta){\Phi}(\diff \Pmu^*)$, where the prior distribution $\Pi(\diff\boldsymbol{\vartheta})$ on the GP parameters  is defined as in formula \eqref{eq:prior_ist} and the prior distribution on the law $\Pmu^*$ is defined as ${\Phi}(\diff \Pmu^*)=\text{DP}(\diff \Pmu^*;\tau)$, namely a Dirichlet process (DP), where $\tau$ is a finite positive measure on Borel sets $B\subset [0,1]^d$ \citep[see e.g.][Ch 4.1]{ghosal2017} {which we hereafter assume absolutely continuous}. 
According to the approximate joint model in \eqref{eq:joint_dist_yx} we have that the posterior distribution for the parameters $(\bfvartheta,\Pmu^*)$ is for all measurable sets $(B\times C)\subset \Theta\times\Psp$ given by
\begin{equation*}\label{eq:posterior_pot_cov}
	\Lambda_n(B\times C)=\Pi_n(B)\Phi_n(C).
\end{equation*}
Since the approximate statistical model arising
from \eqref{eq:joint_dist_yx} postulates independence among the exceedances and the concomitant covariates and given the independence between the prior distributions, then the posterior distribution $\Lambda_n$ splits into the product between the posterior distribution $\Pi_n$ of $\bfvartheta$, given in \eqref{eq:posterior_pot}, and the posterior distribution of $\Pmu^*$, which due to conjugacy property of the DP prior (see Ch 4.6 in \cite{ghosal2017}) becomes $\Phi_n(C)=\text{DP}(C;\tau+k\bbP_n^*)$, that is a Dirichlet process with parameter $\tau+k\bbP_n^*$, where $\bbP_n^*(\cdot)=k^{-1}\sum_{i=1}^k\indic(\bfX_{n-i+1,n}\in\cdot)$ is the empirical measure associated to the covariates concomitant to peaks. For this reason we can {initially} handle the two posterior distributions separately. The inference about $\bfvartheta$ via {$\Pi_n$} is fully described in Section \ref{sec:pot_theory}. Then, we are left here to describe the inference on $\Pmu^*$ via $\Phi_n$, discuss its asymptotic properties and, most importantly, determine the resulting theory for inference on the extreme quantiles of the conditional distribution. This is done by explicitly accounting for the fact that  $\Pi_n$ and $\Phi_n$ are dependent random measures which, however, become increasingly close to Gaussian measures with asymptotically independent random means and deterministic variance as $n \to \infty$  (see Corollary 3.18, Remark 3.19 in the supplement).

The asymptotic theory from the available Bayesian non-parametric literature  (see Ch. 6--12 in \cite{ghosal2017}) cannot be directly applied to the posterior distribution $\Phi_n$, although it is a standard Dirichlet process, as $\Pmu^{*k}$ is a misspecified model for the concomitant covariates associated to the peaks.
We establish here the asymptotic theory of $\Phi_n$ under misspecification, provided that some weak conditions are satisfied. 
Estimation results rely on the control of the convergence speed of the first-order condition in \eqref{eq:tail_prop} through the following second-order condition.
\begin{cond}\label{cond:second_order_c}
There is a  nonincreasing $A_1(t)$, such that $A_1(t)\downarrow 0$ as $t\to\infty$ and
%
$$
\sup_{\bfx\in [0,1]^d}  \left|\frac{1-F^{(0)}_{\bfx}(y)}{1-F_0(y)}-c_0(\bfx)\right|=O\left(A_1\left(\frac{1}{1-F_0(y)}\right)\right), \quad y\to y^*.
$$
\end{cond}

A crucial step is the estimation of the scedasis function $c_0(\bfx)$. For this purpose, we use the fact that whenever $c_0$ is positive and continuous we have
$$
c_0(\bfx)=\lim_{n\to\infty}\frac{\Pmu_0^*(B_n)}{\Pmu_0(B_n)},
$$
where $B_n$ is a sequence of sets containing $\bfx$ and with a decreasing volume.  Accordingly, if $\Pmu_0(B_n)$ was known, the Dirichlet Process prior on $\Pmu^*$ would induce a prior on $\Pmu^*(B_n)/\Pmu_0(B_n)$ and its posterior could be used to infer the scedasis function at $\bfx$. However, $\Pmu_0(B_n)$ is unknown in practice. Then, in the sequel we regard it as a nuisance parameter and assess it by the estimator $\widehat{p}_n\equiv \widehat{p}_n(\bfx) =n^{-1}\sum_{i=1}^{n}\indic(\bfX_i\in B_n)$. We obtain then a data dependent prior on $c(\bfx):=\Pmu^*(B_n)/\widehat{p}_n$, for a given $\bfx$, and we establish the asymptotic properties of its posterior distribution $\Psi_n$ to guarantee a high accuracy of Bayesian inference on $c_0(\bfx)$.

We {first} focus on the situation where $B_n=B(\bfx, r_n)=\{\bfy\in[0,1]^d:\|\bfy-\bfx\|\leq r_n\}$ is the ball of center $\bfx$ and radius $r_n$ and we consider two possible ways of selecting the radius $r_n$: the deterministic one where ball volume is the same for all $\bfx\in[0,1]^d$, we call the resulting estimation procedure {\it kernel} based method with bandwidth ${bw}=r_n$; the data-dependent one where  the ball volume changes depending on $\bfx$ in order to estimate $c(\bfx)$ with a fixed number {$K$} of surrounding points, we call the resulting estimation procedure {\it $K$-nearest neighbours} (KNN) based method.
 Their definition is made precise in the next result.
\begin{theorem}\label{theo:c(x)}
Assume $\Pmu_0$ is absolutely continuous with a density $p_0$ which is positive and locally 
Lipschitz continuous  on {$(0,1)^d$}. Let $c_0$ be a
locally Lipschitz continuous function.
Let $B_n=B(\bfx,r_n)$, where $\bfx \in (0,1)^d$ and
$$
r_n=R\left(\frac{K}{n}\right)^{1/d}\left(1+o\left(\sqrt{\frac{n}{Kk}}\right)\right)\; \mbox{or} \; r_n=\min\left\{h>0\colon \sum_{i=1}^{n}\mathds{1}(X_i\in B(\bfx,h))\geq K\right\},
$$
with $R>0$. Assume $k=o(n)$, $K=o(n)$, $n=o(kK)$ 
and $(K/n)^{1/d}(kK/n)^{1/2}=o(1)$. Assume also that Condition \ref{cond:second_order_c} is satisfied and $kA_1(n/k)\to0$ as $n\to\infty$. 
\begin{itemize}
	\item (Contraction rate) 
	{T}here is a $R'>0$ such that for all sequences $\varepsilon_n$, satisfying $\varepsilon_n\to 0$ and  $(kK/n)^{1/2}\varepsilon_n\to\infty$ as $n\to\infty$, we have that
	\begin{equation*}\label{eq:concentration-rate-c(x)}
		\Psi_n\left(\left\{c(\bfx):\left|c(\bfx)-c_0({\bfx})\right|>\varepsilon_k\right\}\right)= O_\bbP\left(e^{-R'(kK/n)\varepsilon_n^2}\right).
	\end{equation*}
	\item (Bernstein-von Mises) 
	We have 
	\begin{equation*}\label{eq:BvM-c(x)}
		\sup_{B\subset {\mathbb{R}}}\left|\Psi_n\left(\left\{c(\bfx):v_n\left(c(\bfx)-c_0(x)\right)\in B\right\}\right)-\mathcal{N}(B; \Delta_n, c_0(\bfx))\right|=o_\bbP(1),
	\end{equation*}
	where $B$ is any Borel subset of ${\mathbb{R}}$, $v_n=\sqrt{R^{''}kK/n}$, $\Delta_n=v_n(\widehat{p}_n^*/\widehat{p}_n-c_0(\bfx))\stackrel{d}{\to}\mathcal{N}(0,c_0(\bfx))$, $\widehat{p}^*_n{=}k^{-1}\sum_{i=1}^k\mathds{1}(\bfX_{n-k+i,n}\in B_n)$ and {$R^{''}=R^d\pi^{d/2}p_0(\bfx)/\Gamma(1+d/2)$} and $R^{''}=1$ with the above left-hand and right-hand side definition of $r_n$, respectively. 
	\item (Coverage probability) For any 
	$\alpha\in(0,1)$ we have
	$$
	\bbP\left(\left\{c_0(\bfx)\in (\Psi_{n}^{\leftarrow}(\alpha/2), \Psi_{n}^{\leftarrow}(1-\alpha/2))\right\}\right)=1-\alpha+o(1),
	$$
	where, for $a\in(0,1)$, $\Psi_{n}^{\leftarrow}(1-a)$ is the $(1-a)$-quantile of $\Psi_n$.
\end{itemize}
\end{theorem}
Similarly to the unconditional case, assessing the risk associated to extreme events of a certain phenomenon that takes place when other concomitant dynamics reach certain levels is of primary interest in practical problems. This task can be achieved computing extreme conditional quantiles, namely the quantiles of $F_{\bfx}$ corresponding to a small exceedance probability   $p=o(k/n)$,
which can be approximated for large $n$ by the right-hand side of formula \eqref{eq:conditional_quantile}. 
That expression seen as a function of $\bfvartheta_0$ and $c_0(\bfx)$, is a continuous map $T_n: \Theta\times\Real\to\Real$, for any given $n$, $k$, $p$ and $\bfx$,  and therefore {$\Lambda_n$ induces a posterior distribution $\Lambda_n\circ T_n^{-1}$} on the approximate extreme conditional quantile. 
Since our procedure is based on the frequentist estimation of $U_0(n/k)$ via $Y_{n-k,n}$, next we are going to denote the map, the posterior distribution and the extreme conditional quantile as $\widetilde{T}_n$,  ${\widetilde{\Lambda}_n}:={\Lambda}_n\circ \widetilde{T}_n^{-1}$ and 
\begin{equation}\label{eq:cond_extr_quantile}
Q_{\bfx}(p)=Y_{n-k,n}+ H^{\leftarrow}_{\bfvartheta}(1-np/(kc(\bfx))).
\end{equation}
\begin{corollary}\label{theo:quantile_posterior_2}
Assume that conditions of Theorems \ref{theo:main_posterior} and \ref{theo:c(x)} are satisfied and, in the special subcase where $\rho=0$, further  assume that $\gamma_0<0$.  
{Let $p=o(k/n)$ be such that 
$\log(k/np)=o(\sqrt{k})$ and
$$
q_{\gamma_0}(k/np) (k/np)^{-\gamma_0}\sqrt{K/n}\to \omega \in [0,\infty], \quad n\to\infty.
$$
If $\omega=\infty$, define $v_n$ as in Corollary \ref{theo:quantile_posterior}, otherwise set 
$$
v_n=(c_0(\boldsymbol{x}))^{1-\gamma_0}\frac{\sqrt{R''kK/n}}{(k/np)^{\gamma_0}a_0(n/k)}.
$$
Then, for all $\bfx\in(0,1)^d$}:
\begin{itemize}
	\item (Contraction rate) 
	{T}here is a $R^{'''}>0$ such that, for all sequences $\varepsilon_n\to0$ satisfying  $\sqrt{kK/n}\varepsilon_n\to \infty$ and $\sqrt{K/n}\varepsilon_n\log(k/np)\to 0$ as $n\to\infty$, we have 
	$$
	{\widetilde{\Lambda}_n}\left(\left\{Q_{\bfx}(p)\in \Real: \left|\frac{Q_{\bfx}(p)-F^{(0)\leftarrow}_{\bfx}(1-p)}{
	{a_0(n/k)q_{\gamma_0}(c_0(\bfx)k/(np))}
	}\right|>{\widetilde{\varepsilon}_n}\right\}\right)= O_\bbP\left(\exp\left(-R^{'''} {(kK/n)}
	\varepsilon_n^2\right)\right)
	$$
	where $\widetilde{\varepsilon}_n=\varepsilon_n (np/k)^{-\gamma_0}  /q_{\gamma_0}(c_0(\bfx) k/(np))$ if $\omega<\infty$, while $\widetilde{\varepsilon}_n= \varepsilon_n\sqrt{K/n} $ if $=\infty$.
	\item (Bernstein-von Mises) 
	{W}e have
	$$
	\sup_{B\subset \Real}\left|{\widetilde{\Lambda}_n}\left(\left\{Q_{\bfx}(p)\in \Real:v_n\left(Q_{\bfx}(p)-F^{(0)\leftarrow}_{\bfx}(1-p)\right)\in B\right\}\right)-\mathcal{N}(B;{\widetilde{\Xi}_n, \widetilde\Omega})\right|=o_\bbP(1),
	$$
	where $B$ is any Borel subset of $\Real$, 
	{ see Section 3.8 of the supplement for  $\widetilde{\Xi}_n$ and $\widetilde{\Omega}$}.
	\item (Coverage probability) For any  
	$\alpha\in(0,1)$, if $\lambda=0$ we have
	$$
	\bbP\left(\left\{F_{\bfx}^{(0) \leftarrow}(1-p)\in ({\widetilde{\Lambda}_{n}}^{\leftarrow}(\alpha/2), {\widetilde{\Lambda}_{n}}^{\leftarrow}(1-\alpha/2))\right\}\right) = 1-\alpha+o(1),
	$$
	%
	where, for $a\in(0,1)$, $\widetilde{\Lambda}_{n}^{\leftarrow}(1-a)$ is the $(1-\alpha)$-quantile of {$\widetilde{\Lambda}_n$}.
\end{itemize}
\end{corollary}
We complete this section discussing the more complicated case of functional estimation{. Let $\mathcal{G}_n \subset (0,1)^d$ be a grid of  $N \equiv N(n) \to \infty$ points such that $(B(\bfx,r_n), \bfx \in \mathcal{G}_n)$ covers $[0,1]^d$, where $r_n$ is deterministic with $R>0$. 
	Define the piecewise constant function $c(\bfx')=\Pmu^*(B(\bfx,r_n))/ {\widehat{p}_n(\bfx)}	$ for $\bfx =  \mathop{\mathrm{argmin}}_{\bfx \in \mathcal{G}_n}\Vert \bfx' - \bfx \Vert_\infty$, with the convention $c_0(\boldsymbol{x}')=0$ if $
\widehat{p}_n(\bfx)=0$. For simplicity, we denote its posterior distribution still by $\Psi_n$.
\begin{theorem}\label{th:func}
Work under Theorem \ref{eq:concentration-rate-c(x)} conditions, with deterministic radius $r_n$ and Lipschitz continuous functions $p_0, c_0$ on $[0,1]^d$. Set $\delta_0=\inf_{\boldsymbol{x}'\in [0,1]^d}p_0(\boldsymbol{x}')$. 
Then, there is a $R''''>0$ such that, for all positive  $\delta_n=\delta_0 +o(1), \, \varepsilon_n=o(1)$ satisfying 
		 $r_n=o(\delta_n)$, $\varepsilon_n \sqrt{\delta_n kK/(n \log N) } \to \infty$ 
		as $n \to \infty$,  the posterior satisfies
		\[
		\Psi_n \left(\left\lbrace
		c: \,  \Vert c-c_0\Vert_{\delta_n} > \varepsilon_n
		\right\rbrace
		\right) =O_{\mathbb{P}}\left( \exp\left(
		 - R'''' \delta_n (kK/n) \varepsilon_n^2 
		\right) \right)
		\]
		where $\Vert c-c_0\Vert_{\delta_n}= \sup_{\bfx' : \, p_0(\bfx')>\delta_n }|c(\bfx')-c_0(\bfx')|$.
\end{theorem}	
Theorem \ref{th:func} provides uniform contraction rates of the posterior distribution of $c(\bfx)$ over a collection of covariates values $\bfx$, for which the density $p_0(\bfx)$ is positive, {and is therefore a stronger result than 
Theorem	 \ref{theo:c(x)}}. As soon as $p_0$ is bounded away from $0$, our findings cover the $L^\infty$-functional contractions rates, ensuring high accuracy of $\Psi_n$-based inference on the scedasis function.

We finally discuss inference on the distribution function. For any $\bfx\in[0,1]^d$, let $P(\bfx)$ and $P^*(\bfx)$ be the cumulative distribution functions of $\Pmu$ and $\Pmu^*$, respectively, and $P_0(\bfx)$ and $P_0^*(\bfx)$ true counterparts. We {accodingly set} $\mathbb{P}_n^{\circ}(\bfx)=n^{-1}\sum_{i=1}^n\indic(\bfX_i\leq \bfx)$ and $\bbP_n^*(\bfx)=k^{-1}\sum_{i=1}^k\indic(\bfX_{n-k+i,n} \leq \bfx)$. Finally, let $\mathcal{B}(\bfx)$ be a $P_0^*$-Brownian bridge, i.e. a zero-mean Gaussian process with covariance
$$
\Expect[\mathcal{B}(\bfx_1)\mathcal{B}(\bfx_2)]=P_0^*(\min(\bfx_1,\bfx_2))-P_0^*(\bfx_1)P_0^*(\bfx_2),\quad \bfx_1,\bfx_2\in[0,1]^d
$$
and $\widetilde{\Phi}$ be its probability law.} The following result establishes that the posterior distribution $\widetilde{\Phi}_n$ on $P^*$ converges to $\widetilde{\Phi}$ in the functional sense and is therefore also asymptotically Gaussian. 
Note that $\widetilde{\Phi}_n$ and $\widetilde{\Phi}$ are Borel probability measures on the complete and separable Skorohod space $D([0,1]^d, d_0)$, defined in \cite{neuhaus71}. We denote with $\nu(\cdot; \cdot)$ the L\'evy-Prohorov metric \citep[][p. 488]{ghosal2017}, metrising weak convergence over separable spaces.
\begin{theorem}\label{prop:BvM-cdf}(Berstein-von Mises) { Work under Condition \ref{cond:second_order_c}.  Let $\Pmu_0$, $\tau$ be absolutely continuous, with $\tau$ having positive density over $[0,1]^d$. Then, if $k A_1(n/k)=o(1)$ and $k=o(n)$ as $n \to \infty$, it holds that}
$$
\nu \left(
\widetilde{\Phi}_n \left(\{
P^*: \sqrt{k}(P^*-\mathbb{P}_n^*)
 \in \, \cdot \,\} \right); \widetilde{\Phi}
\right)=o_{\bbP}(1).
$$
As a consequence, as $n \to \infty$
$$
\nu\left(
\widetilde{\Phi}_n \left( \{
P^*: \, \Vert \sqrt{k}(P^*-\mathbb{P}_n^*)
\Vert_\infty \in \, \cdot \, \} \right) \, ; \,
\widetilde{\Phi}\left(\{\widetilde{P}: \,  \Vert \widetilde{P} \Vert_\infty \in \, \cdot \, \}  \right)
\right)=o_{\bbP}(1).
$$
\end{theorem}
The practical utility of Theorem \ref{prop:BvM-cdf} is for instance to enable the derivation of a test statistic 
to verify whether the concomitant covariates have a significant effect on the extremes of the response variable. When the covariates have no effect on the extremes of the response variable we have that $c_0(\bfx)=1$ for all $\bfx\in[0,1]^d$. On this basis we consider then the system of hypotheses
$$
\mathscr{H}_0: P_0^*= P_0 \quad \mbox{versus} \quad \mathscr{H}_1: P_0^*\neq P_0,
$$
where again $P_0$ and $P_0^*$ are the cumulative distribution functions of $\Pmu_0$ and $\Pmu_0^*$, respectively. For testing the validity of $\mathscr{H}_0$ we consider a Kolmogorov-Smirnov-type of test as in \cite{einmahl2016} and we rely on the Berstein-von Mises result in Theorem \ref{prop:BvM-cdf}  to draw many samples from the posterior distribution and assess then the critical value. Specifically, in order to perform our hypothesis testing we consider following scheme:
\begin{enumerate}
\item Compute the test statistic $\mathbb{S}=\sqrt{k}\|\mathbb{P}_n^*-\mathbb{P}_n^{\circ}\|_{\infty}$;
\item Draw independent samples $(\mathcal{P}^*_m)_{1\leq m\leq M}$ from the posterior distribution {$\text{DP}(\tau+k\bbP_n^*)$}, for a large value $M$ and then compute the corresponding distributions $(P^*_m)_{1\leq m\leq M}$ and statistics
$$
\mathbb{S}_m=\sqrt{k}\| P_m^*-\mathbb{P}_n^*\|_\infty,\quad m=1,\ldots M;
$$
\item For any $\alpha\in (0,1)$, compute the $(1-\alpha)$-quantile of $(\mathbb{S}_m)_{1\leq m\leq M}$ denoted by $\widehat{Q}_\mathbb{S}(1-\alpha)$. Finally, reject $\mathscr{H}_0$ if $\mathbb{S}>\widehat{Q}_\mathbb{S}(1-\alpha)$.
\end{enumerate}
Thanks to Theorem~\ref{prop:BvM-cdf} we have that asymptotically the significance level of such a hypothesis test is $\alpha$, as $n\to\infty$ and for $M\to\infty$ . We remark that the empirical quantile $\widehat{Q}_\mathbb{S}(1-\alpha)$ is an estimate of the $(1-\alpha)$-quantile of the sup-norm $\|\mathcal{B}\|_\infty$ {of} the $P_0^*$-Brownian bridge $\mathcal{B}$. 

%
%
\subsection{Asymptotic theory of predictive distribution}\label{sec:predictive_EQ}
%
%

One of the most prominent statistical problems is the prediction of certain events in a regression framework. Accurate predictions of  severer extreme events than those occurred in the past is far from being trivial. Whenever this is possible the resulting benefit is huge, because of
the strong impact on real life that such events have. Here we want to go beyond the inference obtained by the posterior distribution {$\widetilde{\Lambda}_n$} on an extreme conditional quantile $Q_{\bfx}(p)$, for a very small $p$, see Section \ref{sec:posterior_EQ}.

This section aim is to propose a probabilistic forecasting method in an extreme regression type of framework. 
Given a past sample $(\bfY_n, \bfX^{(n)})$, where $ \bfX^{(n)}=(\bfX_1,\ldots,\bfX_n)$, let $(Y^\star, \bfX^\star)$ be an independent out-of-sample response variable and covariate vector, representative of future events. We consider the conditional distribution $F^\star_{0,n}(y\mid\bfx)=\bbP(Y^\star\leq y \mid Y^\star>U^{(0)}_{\bfx}(n/k), \bfX^\star=\bfx, \bfX_n, \bfY_n)$, for all $y>U^{(0)}_\bfx(n/k)$ and $\bfx\in(0,1)^d$.
 Similar to Section \ref{sec:predictive_POT}, we leverage the results from Section \ref{sec:posterior_EQ} (see Remark 4.3 in the supplement) in a regression setting to perform forecasting using the posterior predictive distribution.
A possible estimator of  $F_{0,n}^{\star}(\cdot\mid\bfx)$ is given by
\begin{equation}\label{eq:pred_covariates}
\widehat{F}_n^{\star}(y\mid\bfx)={\int_{\Theta\times\mathscr{P}}}
H_{\gamma}
\left(
\frac{y - Y_{n-k,n}}{\sigma\,(c(\bfx))^\gamma}-\frac{1-(c(\bfx))^{-\gamma}}{\gamma}
\right)
\Phi_n(\diff c(\bfx)){\Pi_n(}\diff \bfvartheta),
\end{equation}
see also Remark 4.3 in the supplement for further technical details  justifying its construction.
We recall again that for any $v\geq 1$ the Wasserstein distance of order $v$ satisfies the scaling property
$$
W_v(\widehat{F}_n^{\star}, F_{0,n}^{\star})=a_0(n/k)W_v(\widehat{F}_n^{\star}(a_0(n/k)\cdot\mid\bfx),F_{0,n}^{\star}(a_0(n/k)\cdot\mid\bfx)),
$$
where on the left-hand side we omit conditioning on $\boldsymbol{x}$ for brevity.
Next result {establishes} the Wasserstein consistency of the predictive distribution. This is important as it guarantees that predictions based on the posterior predictive distribution $\widehat{F}_n^{\star}(\cdot\mid\bfx)$ are increasingly accurate for increasing sample size{, whatever is the value of $\bfx$ on $(0,1)^d$}.
\begin{theorem}\label{theo:pred_wass_cov}
Assume that the conditions of Theorem \ref{theo:pred_wass} and \ref{theo:c(x)} are satisfied.
Then, for all sequences $\varepsilon_n\to0$, satisfying  $(kK/n)\varepsilon_n^2/\log(kK/n)\to \infty$ as $n\to\infty$, we have
$$
\frac{W_v(\widehat{F}_n^{\star}, F_{0,n}^{\star})}{a_0(n/k)}=O_\bbP(\varepsilon_n).
$$
\end{theorem}
%

%
%
\section{Simulation experiments}\label{sec:simulations}
We assess the finite sample performance of Bayesian inference based on the posterior distributions introduced in Sections \ref{sec:pot_theory} and \ref{sec:ext_cond_Q}. First, we provide a brief overview of the computational methods used to sample from these posteriors. Given that unconditional analysis is typically of narrower scope than conditional analysis in applications, we then summarize the key findings for the posterior $\Pi_n$ of the parameter $\bfvartheta$ and $\widetilde{\Pi}_n$ of the extreme quantile $Q(p)$ for brevity, with a full description available in Section 6.1 of the supplement. Finally, we present a detailed analysis of the posterior distributions $\Psi_n$ for the scedasis function $c(\bfx)$ and {$\widetilde{\Lambda}_n$} for the extreme conditional quantile $Q{\bfx}(p)$.
\subsection{Posterior distribution computation}\label{sec:computation}
The analytical expression of the posterior $\Pi_n$  is unknown in closed-form. Sampling from it is however viable using MCMC computational methods. The adaptive random-walk Metropolis-Hastings algorithm \citep{haario2001} and its Gaussian random-walk version with Robbins--Monro process optimal scaling (see \cite{garthwaite2016}), is readily implementable and computationally efficient. {It has already been successfully exploited
by \cite{padoan2022empirical}, with the block maxima approach, where extensive simulation experiments demonstrate that an accurate {inference}  is achievable {through the} posterior distribution{,} which complies with the corresponding theoretical findings.}
To save space we provide the full description of such sampling procedure in Section 5 of the supplement. The sampling $Q\sim\widetilde{\Pi}_n$ is achieved as a by product of first sampling $\bfvartheta\sim\Pi_n$ and exploiting the transformation  $Q(p)=Y_{n-k,n}+H^{\leftarrow}_{\bfvartheta}(1-np/k)$, for a small $p\in(0,1)$. Let $\bfvartheta_1^*,\ldots,\bfvartheta_N^*$ be a sample from $\Pi_n$, then a Monte Carlo approximation of the posterior predictive distribution in \eqref{eq:predictive_dist} is 
$$
\widehat{F}^{\star}_n(y)\approx \frac{1}{N} \sum_{i=1}^N H_{\bfvartheta_i^*}(y-Y_{n-k,n})
$$
Since  $\bfvartheta$ and $\Pmu^*$ are independent with distribution $\Pi_n$ and $\Phi_n$, see Sections \ref{sec:proportional_tail} and \ref{sec:posterior_EQ} for details, and given that $\Phi_n$ is a Dirichlet process, the computation of the density and other related quantities of the Dirichlet-multinomial distribution or sampling from it, is readily done using the {\tt R} package {\tt extraDistr} \cite{Tymoteusz20}.
 	
{Conditionally to the data sample, $\bfvartheta$ and $c(\bfx)$ are independent}, then sampling $Q_{\bfx}\sim\widetilde{\Phi}_n$ is achieved sampling first $\bfvartheta\sim\Pi_n$ and independently $c(\bfx)\sim\Psi_n$, for any given $\bfx\in[0,1]^d$, and then transforming them by the formula \eqref{eq:cond_extr_quantile}. Finally, let $\bfvartheta_1^*,\ldots,\bfvartheta_N^*$ be a sample from $\Pi_n$ and $c_1^*(\bfx),\ldots,c_N^*(\bfx)$ be a sample from $\Psi_n$, 
then an approximation of the posterior predictive distribution in \eqref{eq:pred_covariates} is obtained as
$$
\widehat{F}^{\star}_n(y|\bfx)\approx \frac{1}{N} \sum_{i=1}^N H_{\gamma_i^*}
\left(
\frac{y - Y_{n-k,n}}{\sigma_{i}^*\,(c_{i}^*(\bfx))^{\gamma_i^*}}-\frac{1-(c_{i}^*(\bfx))^{-\gamma_{i}^*}}{\gamma_{i}^*}
\right), \quad \bfx\in[0,1]^d.
$$
\subsection{Unconditional POT setting}\label{sec:iid_case}
We investigate the behaviour of the posterior distributions $\Pi_n$ and $\widetilde{\Pi}_n$ and the performance of the resulting inference. The investigation relies on a simulation experiment involving nine distributions, three for each domain of attraction: Fr\'echet, Pareto and Half-Cauchy in the Fr\'echet one, Exponential, Gumbel and Gamma in the Gumbel one and Beta, Weibull and Power-law in the Weibull one. To save space, we refer to Section 6.1 of the supplement for a complete description of: the simulation setup, the computational aspects and the collected results. Here is a summary of our findings. Firstly, we study the concentration properties of the posterior distribution $\Pi_n$,  theoretically implied by the consistency result in Theorem \ref{theo:main_posterior}. With all the considered distributions, the empirical posterior distribution is already fairly concentrated around the true parameter value with only $k=20$ exceedances from a sample of size $n=155$.  In the Fr\'echet domain of attraction the posteriors are more spread than those obtained with the other two domains. However, increasing the sample size $n=303,699,2146$ and the number of exceedances $k=30, 50, 100$ the posterior distribution shrinks considerably and in the last case concentrates  very much in proximity to the true parameter values. These results support the asymptotic concentration properties in Theorem \ref{theo:main_posterior}.

Secondly, we compute the Monte Carlo coverage probability of symmetric- and asymmetric-$95\%$ credible intervals for $\gamma_0$ and  $a_0(n/k)$ and the extreme quantile $F_0^{\leftarrow}(0.999)$ and of symmetric- and asymmetric-$95\%$ credible regions for $\bfvartheta_0$. Overall, with all the models in the three domains of attraction the performance is very good, with coverage probabilities that are close to the $95\%$ nominal level already with the smallest  intermediate sequence $k=20$ and sample size $n=155$. 
With the increasing of $n$ and $k$ the coverage probabilities get even closer.
In Fr\'echet and Gumbel domains of attraction, the coverage probabilities relative to $\gamma_0$ and $a_0(n/k)$ are almost the same. In the Weibull domain of attraction, the symmetric intervals for $\gamma_0$ are slightly larger than expected. Differently, symmetric intervals for $a_0(n/k)$ {have coverage probability slightly smaller than the nominal level.} In the {Fr\'echet and Gumbel}
domains of attraction, the symmetric intervals for  $F_0^{\leftarrow}(0.999)$ are slightly larger than expected{, while in the Weibull domain a coverage probability above nominal level is expected. Overall, the} asymmetric ones perform much better. In the three domains of attraction the worst coverage probabilities are obtained with the credible region for $\bfvartheta_0$, since a higher-dimensional parameter is more difficult to estimate. However, also in this case the coverage probabilities approach the nominal level as $k$ and $n$ increase. 
Concluding, the valuable theoretical properties are actually verifiable in practice already with moderate sample sizes.
\subsection{Extreme regression setting}\label{sec:iid_ext_regressione}
\begin{figure}[t]
	\centering
	\includegraphics[width = \textwidth/(31)*10]{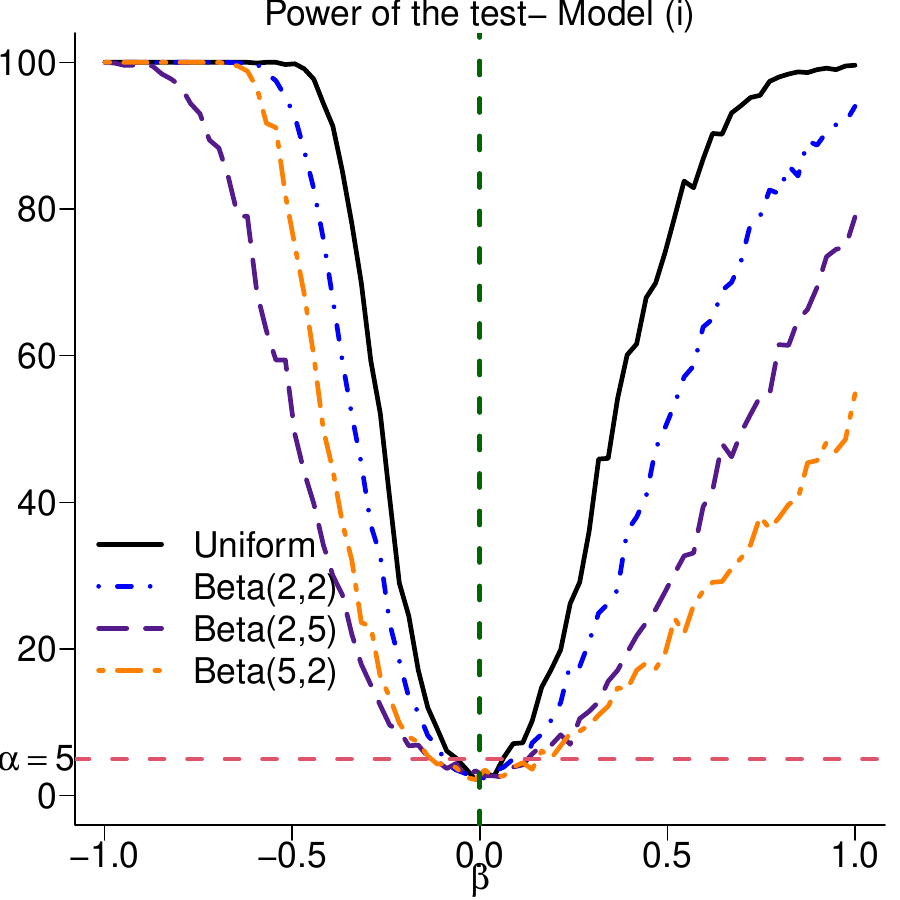}
	\includegraphics[width = \textwidth/(31)*10]{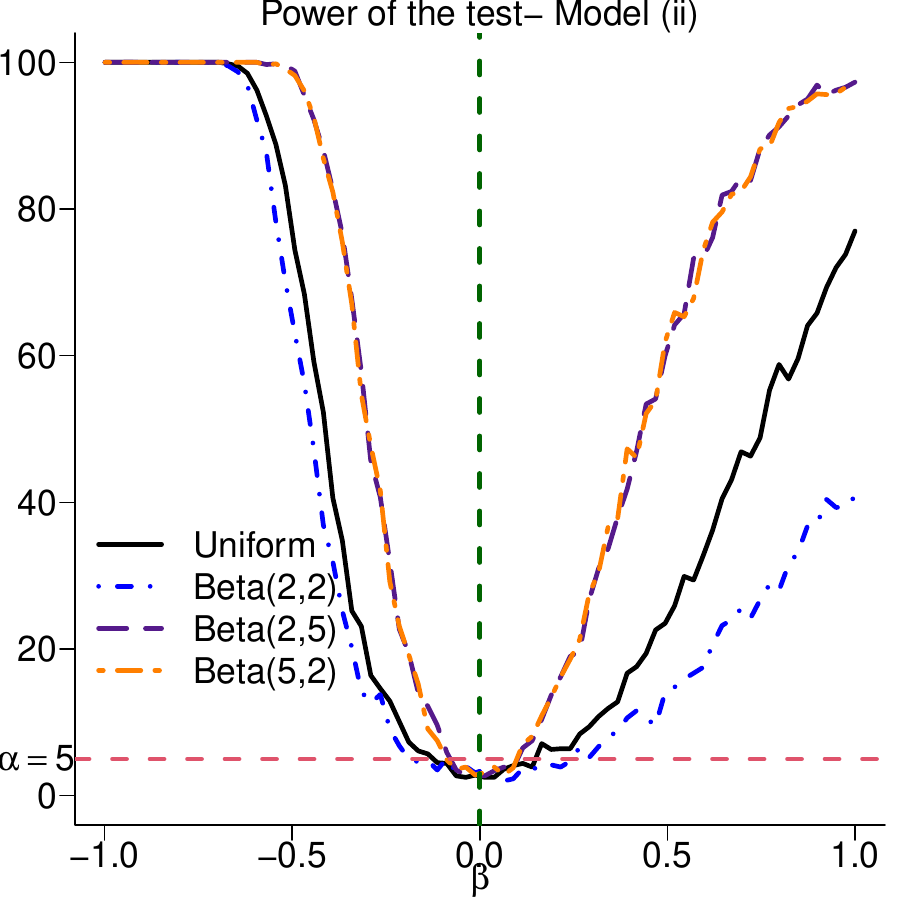}
	\includegraphics[width = \textwidth/(31)*10]{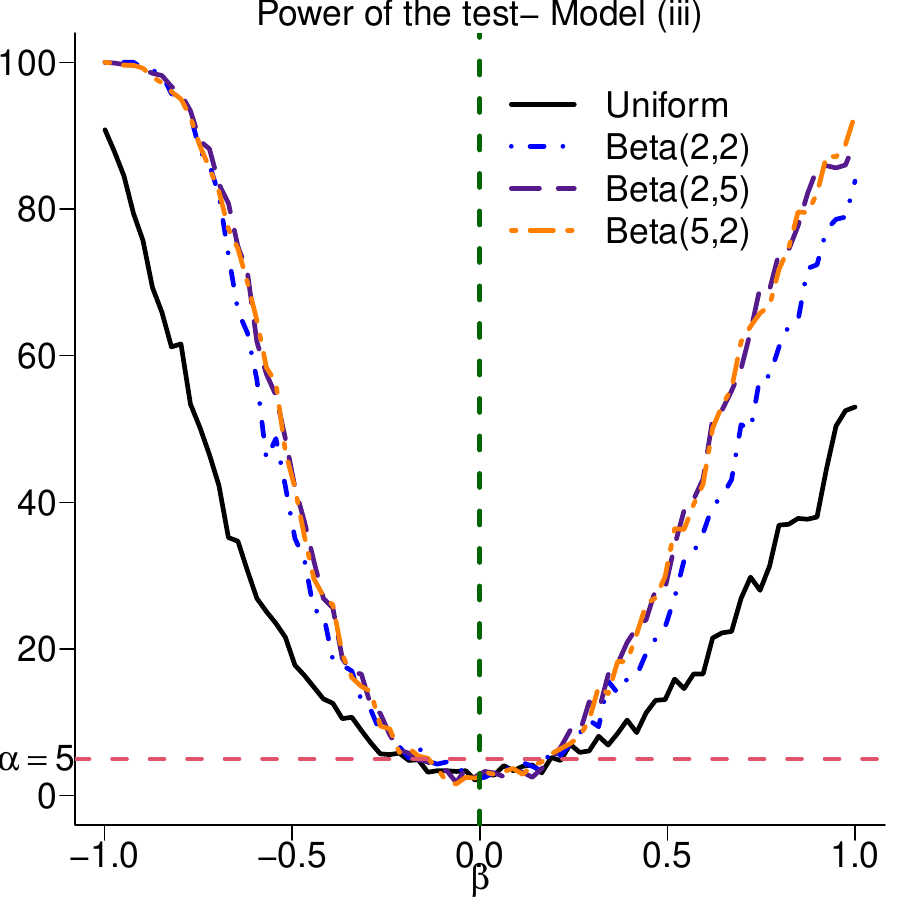}
	\caption{Estimated power functions. Lines report the empirical proportion of simulated samples under {$H_1: P_0^*\neq P_0$} that rejected {$H_0:P_0^* = P_0 $} as a function $\beta$. Dotted red horizontal line is the $5\%$ significance level of the test.}
	\label{fig:power_test}
\end{figure}

We study the behaviour of the posterior distributions $\Phi_n$ and {$\widetilde{\Lambda}_n$} and the performance of the resulting inference. We consider two experiments where the data are generated according the following mechanism. First, we  sample $x_1,\ldots,x_n$ observations from $X_1,\ldots,X_n$ i.i.d. random covariates. To account for the cases that the covariate is scattered over whole [0,1], concentrated on 1/2, concentrated to the left near zero and concentrated to the right near one, we consider the following options for the distribution of $X$: $\mathcal{U}(0,1)$, i.e. uniform on $[0,1]$, Beta$(2,2)$, Beta$(2,5)$ and Beta$(5,2)$, were Beta$(a,b)$ is a Beta distribution with shape parameters $a$ and $b$.
Second, similarly to \cite{einmahl2016}, the $i$th observation $y_i$ is generated from $Y|(X_i=x_i)$ whose distribution is the rescaled Fr\'echet distribution $F_{x_i}(y)=\exp(-c(x_i)/y)$, $y>0$. We consider three possible models for the scedasis function:
\begin{enumerate}
	\item[(i)] [scedasis straight line] $c(x)=(1+\beta\,x)\indic(0\leq x \leq 1)$;
	\item[(ii)] [scedasis broken line] $ c(x)=(1+2\beta\,x)\indic(0\leq x\leq 0.5)$+\\ $(1+2\beta\,(1-x))\indic(0.5<x\leq 1)$;
	\item[(iii)] [scedasis bump function] $c(x)=\indic((0\leq x \leq 0.4)\cup (0.6\leq x \leq 1))+$ $(1+10\beta\,(x-0.4))\indic(0.4 < x \leq 0.5)+$  
	$(1+10\beta\,(0.6-x))\indic(0.5 < x < 0.6)$.
\end{enumerate}
Note that these data generating processes satisfy the proportional tail assumption with scedasis function $c_0$ related to $c$ by $c_0(x)=c(x)/\int_0^1 c(z)f_X(z)\diff z$, where $f_X$ is the density of the covariate.
In the first experiment, $\beta$ is taken to be a sequence of $100$ equally spaced values in $[-1,1]$. For each value of it we simulate $n=5000$ observations, according to the sampling scheme above described, and  we perform the hypothesis testing introduced below Theorem \ref{prop:BvM-cdf}, where we use the  setting: $k=400$, DP prior with parameter $\tau=5\cdot \mathcal{U}(\cdot;0,1)$, where $\mathcal{U}(\cdot;a,b)$ is the uniform measure on $a$, $b$, $M=1000$ independent samples from the DP prior and significance level $\alpha=0.05$. We repeat the sampling and testing steps $N=1000$ times and we compute the rejection rates.
The estimated significance level, as the proportion of simulated samples under {$\mathscr{H}_0:P_0^*(x)=P_0$} that rejects $\mathscr{H}_0$ is: 3.8\% if  $X\sim\mathcal{U}(0,1)$, 4.6\% if $X\sim$Beta$(2,2)$, 3.7\% if $X\sim$Beta$(2,5)$ and 4.1\% if $X\sim$Beta$(5,2)$ with a scedasis straight line; 3.9\% if  $X\sim\mathcal{U}(0,1)$, 3.9\% if $X\sim$Beta$(2,2)$, 4.2\% if $X\sim$Beta$(2,5)$ and 3.7\% if $X\sim$Beta$(5,2)$ with a scedasis broken line; 4.6\% if  $X\sim\mathcal{U}(0,1)$, 4.5\% if $X\sim$Beta$(2,2)$, 4.0\% if $X\sim$Beta$(2,5)$) and 4.2\% if $X\sim$Beta$(5,2)$ with a scedasis bump function.
Figure \ref{fig:power_test} displays the estimated powers of the test, as the proportion of samples simulated under {$\mathscr{H}_1:P_0^*\neq P_0$} that rejects $\mathscr{H}_0$, obtained with different covariate distributions by the black solid line, the blue dotdashed line, the violet dashed line, and the yellow twodashed, respectively, and with the different scedasis models (i)--(iii) from left to the right panel.
Results highlight accurate estimation of $\alpha$ and a good power of the test. The best results are obtained with a scedasis straight line, and in which case $\alpha$ is better estimated with a skewed on the right covariate's distribution. The larger power of test is obtained when covariate is uniformly scattered.
As expected, the  test is less performing with a scedasis broken line and a scedasis bump function, since they are much more complicated functions to estimate. In these cases, $\alpha$ is better estimated if the covariate's distribution is concentrated around $1/2$. While the largest power is obtained if the covariate's distribution is concentrated close the corners $0$ and $1$ with model (ii)  and uniformly scattered instead with model (iii). 

In the second experiment, we simulate $n=5000$ observations from a rescaled Fr\'echet distribution, with $c(x)$ specified as in models (i)--(iii), and we take  $\beta=1,2, 10$ and we apply the sampling procedures of Section \ref{sec:computation} to simulate $N=20,000$ realizations of $\bfvartheta$ and $c(\bfx)$ from $\Pi_n$ and $\Psi_n$  
for $100$ equally spaced values of $x\in[0,1]$. In the first case we use the informative prior on $\gamma$ and the data dependent prior on $\sigma$ described in Section 5 of the supplement and in the second case we use both the kernel based method with bandwidth $bw=0.1$ and the  KNN  based method with $K=750$ neighbours (see Section \ref{sec:posterior_EQ} for details). Again, the DP prior is set with parameter $\tau=5\cdot \mathcal{U}(\cdot;0,1)$. Combining those samples by the transformation \eqref{eq:cond_extr_quantile} we obtain a sample from $\widetilde{\Lambda}_n$.
We repeat these steps $M=1000$ times and compute a Monte Carlo approximation of the root mean integrated relative squared error (RMIRSE), i.e.
$$
\text{RMIRSE}=\left(
\Expect
\left(
\int_0^1
\left(
\frac{\widehat{f}_n(x))}{f_0(x)} -1
\right)^2\diff x
\right)
\right)^{1/2},
$$
where the true function $f_0(x)$ is either $c_0(x)$ or $F^{(0)\leftarrow}_x(0.001)$ and the corresponding estimator $\widehat{f}_n(x)$ is either the scedasis posterior mean $\overline{c}_n(x)$ or the extreme conditional quantile posterior mean $\overline{Q}_{x,n}(p)$.
\begin{table}[t!]
\caption{Monte Carlo approximation of the RMIRSE for the posterior mean estimators $\overline{c}_n$ and $\overline{Q}_{x,n}$. Posterior distributions are computed using the kernel and KNN methods, different scedasis models and different covariate's distribution. The sixth and eighth columns report the relative gain of KNN compared to kernel one.}
\label{tab:perfo_scadasis_equantile}
%
%
%
\begin{tabular}{lccc|cc|cc}
\toprule
&  &  & &\multicolumn{2}{c}{RMIRSE - $\overline{c}_n$} &  \multicolumn{2}{c}{RMIRSE - $\overline{Q}_{x,n}$}\\
%
%
Model  & $X$'s distribution & $n$ & $k$ & kernel  & KNN & kernel  & KNN\\
\midrule
			(i)  & $\mathcal{U}(0,1)$ & $5000$ & $400$ & $0.780$ & $-52.8\%$ &  $3.841$ & $-3.7\%$ \\
			& Beta(2,2) & -- & -- & $1.043$ & $-21.4\%$ & $3.853$ & $-2.8\%$ \\ 
			& Beta(2,5) & -- & -- & $1.334$ & $-0.7\%$ & $4.190$ & $~4.5\%$ \\
			& Beta(5,2) & -- & -- & $2.864$ & $~~22.4\%$ & $5.632$ & $~8.5\%$ \\
			\midrule
			(ii)  & $\mathcal{U}(0,1)$ & $5000$ & $400$ & $1.465$ & $15.4\%$ &  $2.064$ & $-2.0\%$ \\
			& Beta(2,2) & -- & -- & $2.012$ & $17.3\%$ & $2.054$ & $-11.6\%$ \\
			& Beta(2,5) & -- & -- & $3.419$ & $~8.6\%$ & $7.074$ & $~41.1\%$ \\
			& Beta(2,5) & -- & -- & $3.493$ & $~9.6\%$ & $7.820$ & $~39.8\%$ \\
			\midrule
			(iii)  &  $\mathcal{U}(0,1)$ & $5000$ & $400$ & $1.316$ & $~0.9\%$ &  $2.143$ & $~5.1\%$ \\
			& Beta(2,2) & -- & -- & $1.580$ & $20.1\%$ & $2.278$ & $11.8\%$ \\
			& Beta(2,5) & -- & -- & $5.714$ & $41.9\%$ & $5.105$ & $11.8\%$ \\
			& Beta(5,2) & -- & -- & $4.777$ & $~31.7\%$ & $5.303$ & $12.8\%$ \\
			\bottomrule
		\end{tabular}
\end{table}
%
%
Table \ref{tab:perfo_scadasis_equantile} reports the results split according to the scedasis model (vertical sections), the different covariate's distributions (along the rows), the function to be estimated ($c_0(x)$ in the fifth and sixth column and $F^{(0)\leftarrow}_x(0.001)$ in the seventh and eight column) and the estimation based method for prior and  posterior construction (kernel and KNN). The RMIRSE obtained with the kernel based method show greater precision when estimating a scedasis straight line in comparison to the other two cases.
The difference among results is also fairly small when estimating $F^{(0)\leftarrow}_x(0.001)$ but the RMIRSE does not highlight a clear superiority obtained with a specific scedasis form.
%
Regardless of the scedasis form, the posterior mean is more accurate when the covariate is uniformly distributed or symmetrically concentrated around $1/2$ than in the other cases. Overall, results suggest that the posterior mean is an accurate estimator for $c_0(x)$ and $F^{(0)\leftarrow}_x(0.001)$. The sixth and eighth columns of Table \ref {tab:perfo_scadasis_equantile} report the relative gain (in percentage) of using the KNN method in place of the kernel one, i.e. (RMIRSE(kernel)-RMIRSE(KNN))/RMIRSE(kernel) $\cdot$100\%. When estimating $c_0(x)$ ($F^{(0)\leftarrow}_x(0.001)$), apart from three (four) cases the remaining one highlight better performance of the KNN method, with a gain that ranges between 0.9\% (4.5\%) to 41.9\% (41.1\%) and therefore on balance it is preferable.
%
%
\begin{figure}[t!]
	\centering
	\includegraphics[width=.30\textwidth, height=.18\textheight]{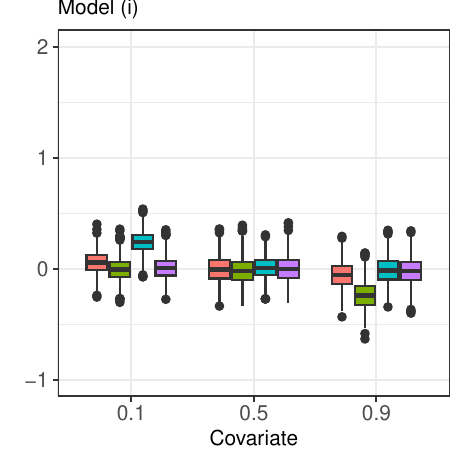}
	\includegraphics[width=.30\textwidth, height=.18\textheight]{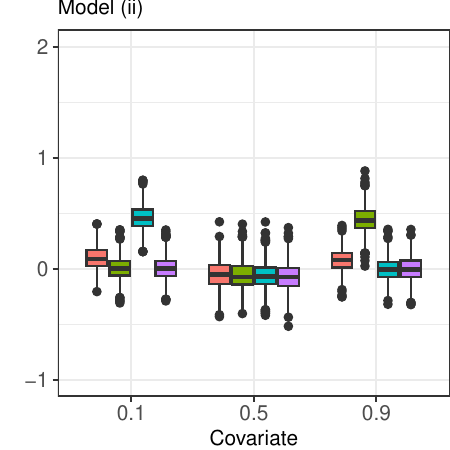}
	\includegraphics[width=.36\textwidth,  height=.18\textheight]{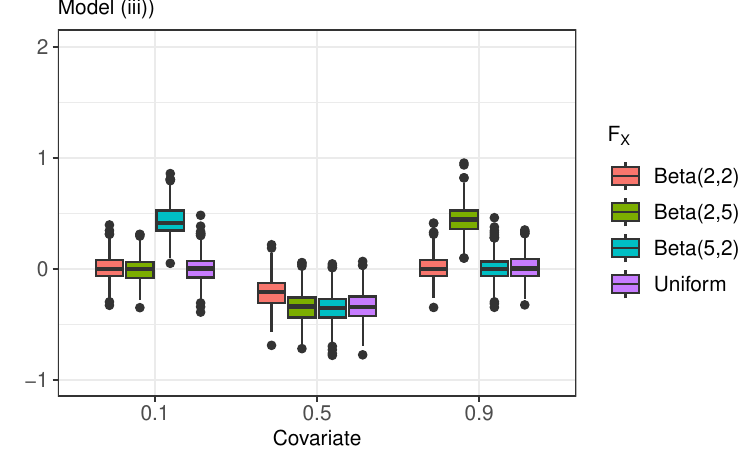}\\
	\includegraphics[width=.30\textwidth, height=.18\textheight]{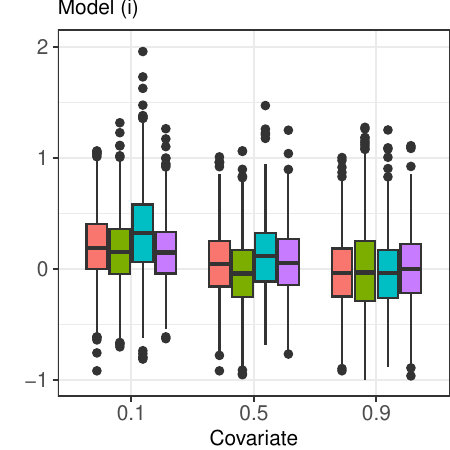}
	\includegraphics[width=.30\textwidth, height=.18\textheight]{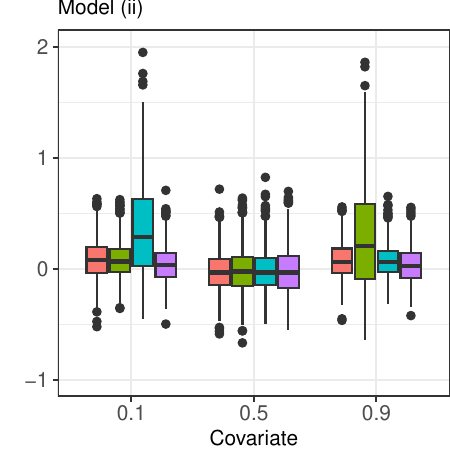}
	\includegraphics[width=.36\textwidth,  height=.18\textheight]{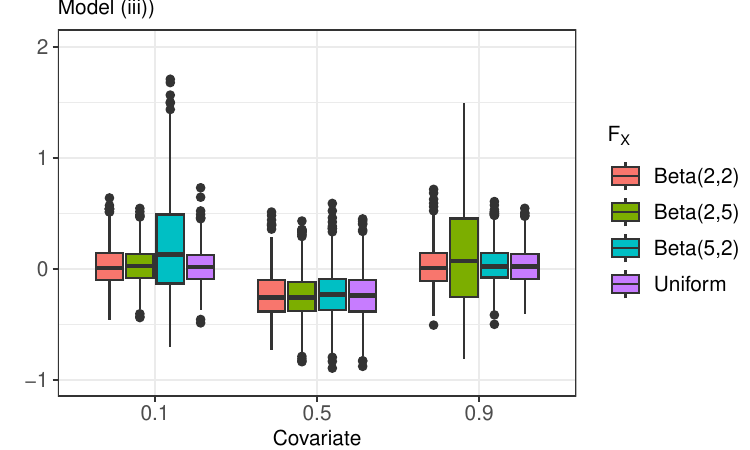}
	\caption{Boxplot of $\overline{c}_n(x)-c_0(x)$ (top panels) and  $\log(\overline{Q}_{x,n})-\log(F^{(0)\leftarrow}_x(0.001))$ (bottom panels), obtained with the scedasis model (i)-(iii) and different covariate's distributions and setting the covariate value  $x=0.1,0.5,0.9$.}
	\label{fig:box_plots}
\end{figure}
%
%

In the sequel we focus on the results obtained with the KNN method. The top panels of Figure \ref{fig:box_plots} display, for three specific values $x=0.1, 0.5, 0.9$, the Monte Carlo distribution of $\overline{c}_n(x)-c_0(x)$, obtained with the $M=1000$ data samples. Results obtained with the model (i)-(iii) are displayed from left to right panel. Each panel reports the results obtained with the different covariate's distributions. Since the boxplots are almost all centred around zero with a small dispersion we can conclude that $\overline{c}_n(x)$ is an accurate estimator of $c_0(x)$. Note that, when the covariate distribution is left (right) skewed, e.g. the Beta(5,2) (Beta(2,5)), $\overline{c}_n(x)$ overestimates a bit $c_0(0.1)$ ($c_0(0.9)$) as expected, since only few data are available around $x=0.1$ ($x=0.9$). With model (iii), $\overline{c}_n$ underestimates a bit $c_0(0.5)$ when the covariate's distribution are not concentrated around $x=1/2$.  
Similar results are obtained when estimating $F^{(0)\leftarrow}_x(0.001)$, see the bottom panels of Figure \ref{fig:box_plots}. The distribution of $\overline{Q}_{x,n}(x)-F^{(0)\leftarrow}_x(0.001)$ is more spread as expected since the estimation of the conditional extreme quantile is harder.

Finally, we compute the Monte Carlo coverage probability of the credible intervals for $c_0(x)$ and $F^{(0)\leftarrow}_x(0.001)$.
To save space, results are reported in Table 4 of Section 6.2 in the supplement. Results are split according to the different scedasis form (vertical sections) and covariate's distributions (along the rows). The column Type indicates with the letter ``A" the coverage of an asymmetric-95\% credible interval (obtained with the quantiles of the posterior distribution) and with the letter ``S" the symmetric version (obtained as $[\widehat{f}(x)_n\pm z_{\alpha/2}\widehat{s}_n(x)]$, where $\widehat{s}_n(x)$ is the posterior standard deviation and $z_{\alpha/2}$ is the standard normal $(1-\alpha/2)$-quantile). With model (i)  the coverages are very close to the  nominal level apart for the case $c_0(0.9)$ ($c_0(0.1)$) when the covariate is Beta(2,5) (Beta(5,2)) distributed, but this is expected as there are few observations close to one (zero). Therefore, a larger sample size is needed to achieve the nominal level all over $[0.1]$. Similar conclusions hold with model (ii) and (iii).

%
%
\section{Application}\label{sec:application}
%
%
%
For comparison, we conduct a similar analysis to that in \cite{einmahl2016} to investigate whether the frequency of financial crises has changed over time. Using the sequence of daily negative log-returns (hereafter, returns) of the Standard \& Poor’s 500 (S\&P 500) index—representing the status of the U.S. financial market—from 1988 to 2012, we focus on the subseries from 1988 to 2007, totaling 5,043 observations. This choice aligns with \cite{einmahl2016}, where the EVI of this shorter series was shown to be time-invariant, in contrast to the full dataset. For simplicity, and following \cite{einmahl2016}, we initially disregard temporal dependence in this analysis.

First we perform the hypothesis test described below Theorem \ref{prop:BvM-cdf}, drawing $M = 1000$ samples from a DP prior with parameter $5\cdot \mathcal{U}(\cdot;0,1)$. The significance level is set to $\alpha = 0.05$, and we choose $k = 210$, a value within the range $[110,250]$, where the EVI estimates remain relatively stable (see Section 7 of the supplement for further discussion). The time coordinate is used as a covariate by mapping trading days to uniformly spaced values in $[0,1]$.
The observed test statistic is $3.604$, while the estimated critical value is $1.276$. Consequently, consistent with the findings of \cite{einmahl2016}, we reject the null hypothesis of a constant scedasis function.
%
%
\begin{figure}[t!]
    \centering
    \includegraphics[width = \textwidth/(31)*10, page=1]{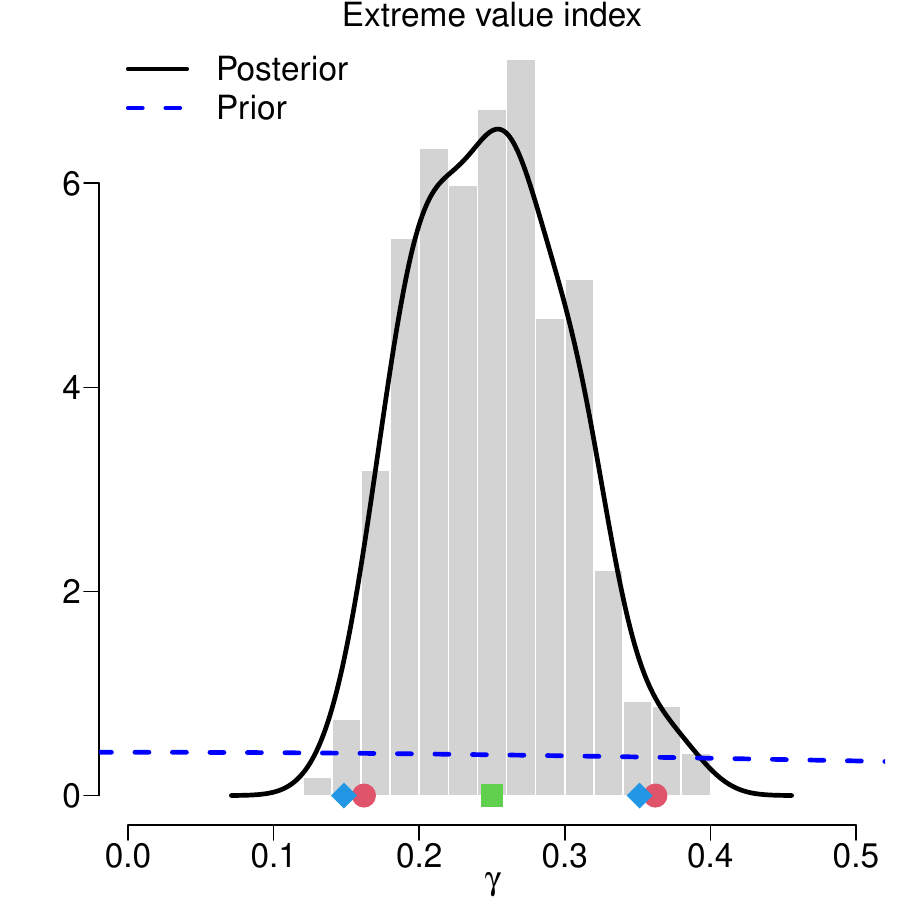}
    \includegraphics[width = \textwidth/(31)*10, page=2]{sp500_bayes_analysis.pdf}
    \includegraphics[width = \textwidth/(31)*10, page=3]{sp500_bayes_analysis.pdf}
    \caption{S\&P 500 return estimation results. EVI Posterior distribution (left panel), estimated scedasis function (middle panel) and loss-returns with conditional extreme quantile estimates and predictive intervals superimposed.}
    \label{fig:S&P500_analysis}
\end{figure}
%
%
%

Next, we compute the empirical versions of $\Pi_n$ and $\Psi_n$ by sampling 20,000 values of $\bfvartheta$ and $c(x)$, with $x\in[0,1]$, using the sampling methods outlined in Section \ref{sec:computation}. Specifically, for $\Pi_n$, we employ a data-dependent prior (see Section 5 of the supplement), while for $\Psi_n$, we utilize the kernel-based method with bandwidth $bw=0.08$ and the KNN-based method with $K=800$. Additionally, the DP prior is set to $5\cdot \mathcal{U}(\cdot;0,1)$.
Figure \ref{fig:S&P500_analysis} presents the results. The left panel shows the posterior distribution of $\gamma$, with  a mean of $0.25$, a standard deviation of $0.052$, and asymmetric (symmetric) 95\% credible intervals of $[0.16,\,0.36 ]$ ($[0.15,\,0.35]$). These results strongly support the assumption of a positive EVI.
In the middle panel, the green solid line represents the posterior mean  $\overline{c}_n(x)$ of the scedasis function, estimated using the kernel-based method (the KNN-based method yields similar results). Consistently with the findings in \cite{einmahl2016}, the posterior mean exhibits a fluctuating pattern, with a pronounced peak around April 2001, followed by a steep increase leading up to 2007.
Compared to the estimation method in \cite{einmahl2016}, our Bayesian approach not only provides an estimate of the scedasis function but also quantifies estimation uncertainty. This is illustrated by the asymmetric (symmetric) 95\% credible intervals, shown as orange dashed and violet dot-dash lines, respectively. The relatively narrow width of these intervals suggests that the period of highest loss risk is statistically likely to have occurred between 2000 and 2002, coinciding with the burst of the dot-com bubble.

Finally, for forecasting purposes, we extend our analysis to a longer time horizon, spanning from 1988 to September 2008. This period includes major financial shocks, such as the significant asset write-downs by major U.S. investment banks in early 2008 and the bankruptcy of Lehman Brothers on September 15, 2008. As before, we transform the trading days into equally spaced values within  $[0,1]$.

The computation of $\Pi_n$ and $\Psi_n$ remains based on returns from 1988 to 2007. However, since $\Psi_n$ is assessed for values of $x$ across the entire interval $[0,1]$, the scedasis function $c$ is now estimated beyond the observed data, covering the full period up to September 2008. Following the approach described in Section \ref{sec:computation}, we approximate the posterior distribution of the extreme conditional quantile, $\widetilde{\Lambda}_n$, with $p=0.001$, as well as the posterior predictive distribution $\widehat{F}^{\star}_n(\cdot \mid x)$ in \eqref{eq:pred_covariates}, using samples drawn from $\Pi_n$ and $\Psi_n$.
The right panel of Figure \ref{fig:S&P500_analysis} illustrates the results. The black solid line represents the returns from 1988 to 2007, while the grey solid line corresponds to returns in the first nine months of 2008. The superimposed green dashed line shows the posterior mean, while the orange dot-dash lines and blue long-dash lines depict the asymmetric 95\% credible from $\widetilde{\Lambda}_n$ and predictive intervals $\widehat{F}^{\star}_n(\cdot \mid x)$, respectively.
Notably, the posterior mean and intervals closely resemble the estimated scedasis function. The credible interval is relatively narrow and, in several instances, does not encompass large losses. In contrast, the predictive interval is significantly wider, capturing most of the major losses observed during the period—particularly the sharp decline following Lehman Brothers' bankruptcy, which falls within the forecasting horizon.

\section*{Acknowledgments}
Simone Padoan is supported by the Bocconi Institute for Data Science and
Analytics (BIDSA) and project MUR - Prin 2022 - Prot. 20227YZ9JK, Italy.

\bibliographystyle{chicago} 
\bibliography{bibliography}
\end{document}